\documentclass[onecolumn]{IEEEtran}

\IEEEoverridecommandlockouts 
\usepackage{cite}
\usepackage{times}  
\usepackage{helvet}  
\usepackage{courier,wrapfig,lipsum}  
\usepackage{url}  
\usepackage{graphicx}  
\usepackage{amsmath,amsfonts,bm,bbm}
\usepackage{graphics, subfigure, times, amssymb, verbatim, amsthm, comment}
\usepackage{color}
\usepackage{hyperref}
\usepackage{algorithm,algorithmic,multirow,hhline}
\usepackage[hang,flushmargin]{footmisc}

\usepackage{subfigure} 
\usepackage{mathrsfs}
\input{mysymbol.sty}

\newcommand{\QED}{\hfill\ensuremath{\blacksquare}}
\renewcommand{\theenumi}{\roman{enumi}}%

\newcommand{\bellman}{{\mathscr{B}}}

\newcommand{\closure}[2][3]{{}\mkern#1mu\overline{\mkern-#1mu#2}}

\newtheorem{assumption}{\hspace{0pt}\bf Assumption}
\newtheorem{lemma}{\hspace{0pt}\bf Lemma}
\newtheorem{theorem}{\hspace{0pt}\bf Theorem}
\newtheorem{proposition}{\hspace{0pt}\bf Proposition}
\newtheorem{corollary}{\hspace{0pt}\bf Corollary}
\newtheorem{remark}{\hspace{0pt}\bf Remark}

\newcommand{\INDSTATE}[1][1]{\STATE\hspace{3mm}}
\newcommand{\INDSTATED}[1][1]{\STATE\hspace{6mm}}

\begin{document}

\title{Policy Evaluation in Continuous MDPs with Efficient Kernelized Gradient Temporal Difference}


\author{
Alec Koppel\textsuperscript{1},
Garrett Warnell\textsuperscript{1},
Ethan Stump\textsuperscript{1},
Peter Stone\textsuperscript{2},
Alejandro Ribeiro\textsuperscript{3}\\
\textsuperscript{1}{U.S. Army Research Laboratory}, \; 
\textsuperscript{2}{The University of Texas at Austin}, \;
\textsuperscript{3}{University of Pennsylvania},\\
{aribeiro@seas.upenn.edu}, 
{\{alec.e.koppel,garrett.a.warnell,ethan.a.stump2\}.civ@mail.mil},
{pstone@cs.utexas.edu}

\author{Alec Koppel, \IEEEmembership{Member, IEEE}, Garrett Warnell, Ethan Stump, Peter Stone, and Alejandro Ribeiro  \IEEEmembership{Member, IEEE}
\thanks{This manuscript was originally submitted for review on December 6, 2017, and supported by the SMART Scholarship, 
ARL DCIST CRA W911NF-17-2-0181, NSF DGE-1321851, Intel DevCloud and Intel Science and
Technology Center for Wireless Autonomous Systems (ISTC-WAS).}
\thanks{\textsuperscript{1} Computational and Information Sciences Directorate, U.S. Army Research Laboratory, 2800 Powder Mill Rd., Adelphi, MD 20783, email:  { \{alec.e.koppel,garrett.a.warnell,ethan.a.stump2\}.civ@mail.mil}.}
\thanks{\textsuperscript{2} Department of Computer Science, University of Texas at Austin, 2317 Speedway, Stop D9500
Austin, Texas 78712-1757 USA, email: {pstone@cs.utexas.edu},  }
\thanks{\textsuperscript{3} Department of Electrical and Systems Engineering, University of Pennsylvania, 200 South 33rd Street, Philadelphia, PA 19104}
}
}
 
\maketitle

\begin{abstract}%
We consider policy evaluation in infinite-horizon discounted Markov decision problems (MDPs) with continuous compact state and action spaces. We reformulate this task as a compositional stochastic program with a function-valued decision variable that belongs to a reproducing kernel Hilbert space (RKHS). We approach this problem via a new functional generalization of stochastic quasi-gradient methods operating in tandem with stochastic sparse subspace projections. The result is an extension of gradient temporal difference learning that yields nonlinearly parameterized value function estimates of the solution to the Bellman evaluation equation.  We call this method Parsimonious Kernel Gradient Temporal Difference (PKGTD) Learning. Our main contribution is a memory-efficient non-parametric stochastic method guaranteed to converge exactly to the Bellman fixed point with probability $1$ with attenuating step-sizes under the hypothesis that it belongs to the RKHS. Further, with constant step-sizes and compression budget, we establish mean convergence to a neighborhood and that the value function estimates have finite complexity. In the Mountain Car domain, we observe faster convergence to lower Bellman error solutions than existing approaches with a fraction of the required memory.
\end{abstract}

\vspace{-4mm}
\section{Markov Decision Processes}\label{sec:mdps}
We consider an autonomous agent acting in an environment defined by a Markov decision process (MDP) \cite{sutton1998reinforcement} with continuous spaces, which is increasingly relevant to emerging technologies such as robotics \cite{kober2013reinforcement}, power systems \cite{scott2014least}, and others. A MDP is a quintuple $(\ccalX, \ccalA, \mathbb{P}, r, \gamma)$, where $\mathbb{P}$ is the action-dependent transition probability of the process: when the agent starts in state $\bbx_t \in \ccalX \subset \reals^p$ at time $t$ and takes an action $\bba_t \in \ccalA$, a transition to next state $\bby_t\in \ccalX$ is distributed according to 
%
$\bby_t \sim \mathbb{P}(\cdot \given \bbx_t, \bba_t).$
%
After transitioning to a particular $\bby_t$, the MDP reveals an instantaneous reward $r(\bbx_t, \bba_t, \bby_t)$, where the reward function is a map $r :\ccalX \times \ccalA \times \ccalX \rightarrow \reals$. 

We focus on \emph{policy evaluation}: control decisions $\bba_t$ are chosen according to a fixed stationary stochastic policy $\pi: \ccalX \rightarrow \rho(\ccalA)$, where $\rho(\ccalA)$ denotes the set of probability distributions over $\ccalA$. Policy evaluation underlies methods that seek optimal policies through repeated evaluation and improvement. In policy evaluation, we seek to compute the \emph{value} of a policy when starting in state $\bbx$, quantified by the discounted expected sum of rewards, or value function $V^{\pi}(\bbx)$:
 \footnote{In MDPs more generally, we choose actions $\{\bba_t \}_{t=1}^\infty$ to maximize the reward accumulation  starting from state $\bbx$, i.e.,
$V(  \bbx, \{\bba_t\}_{t=0}^\infty ) = \mathbb{E}_{\bby}\Big[\sum_{t=0}^\infty \gamma^t r(\bbx_t, \bba_t, \bby_t ) \given \bbx_0 = \bbx, \{\bba_t \}_{t=0}^\infty \Big].$ For a fixed policy $\pi$, the setting of this work, this simplifies to \eqref{eq:value_function}.}\vspace{2mm}
  %
\begin{align}\label{eq:value_function}
\!\!\!V^{\pi}\!(  \bbx ) \!=\! \mathbb{E}_{\bby}\!\Big[\!\sum_{t=0}^\infty \! \!\gamma^t  r(\bbx_t, \bba_t, \bby_t \!)  \!\given \! \bbx_0 \!=\! \bbx, \!\{\!\bba_t \!=\! \pi(\bbx_t)\!\}_{t=0}^\infty \! \Big].
\end{align}
For a single trajectory through the state space $\ccalX$, $\bby_t = \bbx_{t+1}$. The value function \eqref{eq:value_function} is parameterized by a discount factor $\gamma \in (0,1)$ that determines farsightedness. 
%
%
Decomposing the summand in \eqref{eq:value_function} into its first and subsequent terms, and using both the stationarity of the transition probability and the Markov property 
yields the Bellman evaluation equation \cite{Bellman:1957}:
\begin{equation}\label{eq:bellman_value}
V^\pi( \bbx) = \int_{\ccalX} [ r(\bbx,\pi(\bbx),\bby) + \gamma V^{\pi}(\bby)] \mathbb{P}(d\bby \given \bbx, \pi(\bbx))
\end{equation}
for all $\bbx \in \ccalX$.
The right-hand side of \eqref{eq:bellman_value} defines a Bellman evaluation operator $\bellman^\pi : \ccalB(\ccalX) \rightarrow \ccalB(\ccalX )$ over $\ccalB(\ccalX )$, the space of bounded continuous value functions $V:\ccalX\rightarrow \reals$:
%
%
\begin{align}\label{eq:bellman_operator_v}
\!\!\!(\!\bellman^{\pi} V) (\bbx) \! =\!\!    \int_{\ccalX} \!\![ r(\bbx,\pi(\bbx) ,\bby)\! +\! \gamma V^{}(\bby)] \mathbb{P}(d\bby \given \bbx, \pi(\bbx) )
\end{align}
for all $\bbx \in \ccalX$. {Proposition 4.2(b) in \cite{bertsekas1978stochastic}} establishes that the stationary point of \eqref{eq:bellman_operator_v} is $V^{\pi}$, i.e., $( \bellman^{\pi} V^\pi )(\bbx) = V^{\pi}(\bbx)$.
{As a stepping stone to finding optimal policies in infinite MDPs, we seek here to find the fixed point of \eqref{eq:bellman_operator_v}.  Specifically, the goal of this work is stable value function estimation in infinite MDPs, with nonlinear parameterizations that are allowed to be infinite, but are nonetheless memory-efficient.} 

%

{\bf Challenges } To solve \eqref{eq:bellman_operator_v}, fixed point methods, i.e., value iteration ($V_{k+1} = \bellman^{\pi} V_k$), have been proposed \cite{bertsekas1978stochastic}, but {can only be implemented in a memory-affordable manner} when the {value function can be represented by a vector whose length is defined by the number of states and the state space is small enough that the expectation\footnote{\noindent Observe that the integral in \eqref{eq:bellman_value} defines a conditional expectation: $V^\pi( \bbx) = \mathbb{E}_{\bby}[ r(\bbx,\pi(\bbx),\bby) + \gamma V^\pi(\bby)]  \given \bbx, \pi(\bbx)]$.} in $\bellman$ can be computed.} For large spaces, stochastic approximations of value iteration, i.e., temporal difference (TD) learning \cite{sutton1988learning}, circumvent computing expectations. 
Incremental methods (least-squares TD) are an alternative when $V(\bbx)$ is vector-valued \cite{bradtke1996linear}, but extensions to infinite representations require infinite memory \cite{powell2011review}. 

Solving the fixed point problem defined by \eqref{eq:bellman_operator_v} requires surmounting the fact that this expression is defined for {each $\bbx \in \ccalX$}, which for continuous $\ccalX\subset \reals^p$ has \emph{infinitely many} unknowns. This phenomenon is one example of Bellman's curse of dimensionality \cite{Bellman:1957}, and it is frequently sidestepped by parameterizing the value function {using} a finite linear \cite{tsitsiklis1997analysis,melo2008analysis} or nonlinear \cite{bhatnagar2009convergent} basis expansion. Such methods {have} paved the way for the recent success of neural networks in value function-based approaches to MDPs{, but combining TD learning with different parameterizations may cause divergence \cite{Baird95residualalgorithms}: in} general, the representation must be tied to the stochastic update  \cite{jong2007model} to ensure both are stable.

{\bf Contributions } Our main result is a memory-efficient, non-parametric{,} stochastic method that converges to the Bellman fixed point almost surely when it belongs to a reproducing kernel Hilbert space (RKHS). The hypothesis that the value function belongs to a RKHS restricts the relationship between rewards and value to be smooth (see Assumption \ref{as:mean_variance_frechet} [cf. \eqref{eq:td_lipschitz}] ), i.e., large changes in rewards yield large changes in value., {which holds, for instance, when the reward is a potential or navigation function \cite{rimon1992exact,ngicml1999}. Such specifications are known to interact favorably in controller design from a dynamical systems perspective.} We reformulate \eqref{eq:bellman_value} as a compositional stochastic program (Section \ref{subsec:opt}), a topic studied in operations research \cite{shapiro2014lectures} and probability \cite{Korostelev}. These problems motivate stochastic \emph{quasi-gradient} (SQG) methods, i.e.,  two time-scale approaches, to mitigate the fact that the objective's stochastic gradient requires evaluating an expectation \cite{ermoliev1983stochastic}. 

{Two time-scale stochastic methods have a significant history in reinforcement learning in the context of a class of policy learning algorithms called actor-critic \cite{borkar1997actor}, which mix together gradient ascent on the value function \cite{sutton2000policy} with value function estimation through temporal differences \cite{sutton1988learning}. Alternatively, their utility also has been established for policy evaluation in \cite{bhatnagar2009convergent,sutton2009convergent} for finite MDPs or value functions that have finite vector-valued parameterizations. Our work is more closely related to this later context; however, a key point of departure is that we propose to operate directly in a function space, motivated by the fact that policy evaluation in a continuous MDP defines a functional fixed point problem.} Specifically, we use SQG in infinite MDPs. 

In \eqref{eq:bellman_value}, the decision variable is a continuous function, which we address by hypothesizing the Bellman fixed point belongs to a RKHS \cite{kimeldorf1971some}. However, a function in a RKHS has comparable complexity to the {number of training samples processed}, which could be infinite ({an issue} ignored in many kernel methods for MDPs \cite{ormoneit2002kernel,taylor2009kernelized,grunewalder2012modelling,JMLR:v17:13-016,dai2016learning}). {Specifically, it's well known that for a kernelized interpolator defined by a training set of $N$ samples, the complexity is $\ccalO(N)$. Thus, to solve the population problem defined by Bellman's equations, $N\rightarrow \infty$, and thus so is the complexity of the function estimate. We propose to tackle this memory bottleneck by requiring memory efficiency in both the function sample path and in its limit, whose complexity in the worst case is defined by the metric entropy of the state space (Corollary \ref{corollary1}).}

To find a memory-efficient sample path in the function space, we generalize {SQG to RKHSs} (Section \ref{sec:algorithm}), and combine {this generalization} with {greedily-constructed} sparse subspace projections (Section \ref{subsec:polk_projection}). {These subspaces are constructed via matching pursuit \cite{Pati1993,lever2016compressed}, {a procedure} motivated by the {facts} that kernel matrices induced by arbitrary data streams likely violate requirements for convex-relaxation-based sparsity \cite{candes2008restricted}. 
Rather than unsupervised forgetting \cite{engel2003bayes}, we tie the projection-induced error to ensure the stochastic gradient still satisfies a descent property \cite{koppel2019parsimonious}, thus keeping only those dictionary points needed to converge (Section \ref{sec:convergence}). 

{Whereas in \cite{koppel2019parsimonious}, compressed kernel methods are analyzed for supervised learning and a tunable tradeoff between memory and sub-optimality is provided, here we study compressed kernel methods in the context of policy evaluation in reinforcement learning. This later context requires surmounting the technical challenges associated with nested expectations, specifically, that SQG is defined by coupled supermartingales, rather than standard stochastic descent arguments in \cite{koppel2019parsimonious}, and hence has fundamental qualitative and quantitative departures from existing works on RKHS learning.}

We note that this hard-thresholding projection could be applied to other stochastic algorithms in RKHS for reinforcement learning such as \cite{ormoneit2002kernel,taylor2009kernelized,grunewalder2012modelling,dai2016learning}, but applying them to incremental methods (LSTD) \cite{powell2011review,JMLR:v17:13-016} remains elusive since relating the per-step bias caused by sparsification to ensure valid descent directions is elusive.}

{As a result,} we conduct functional SQG descent via sparse projections of the SQG{. This maintains a moderate-complexity} sample path exactly towards $V^*$, which may be made arbitrarily close to {the Bellman fixed point} by decreasing the regularizer.
By generalizing the relationship between SQG and supermartingales in \cite{wang2017stochastic} to Hilbert spaces, we establish that the sparse projected SQG sequence converges almost surely to the Bellman fixed point with decreasing learning rates, and converges in mean while maintaining finite complexity when constant learning rates are used (Section \ref{sec:convergence}). We then empirically evaluate the proposed value function approximation method on the discrete Mountain Car domain in Section \ref{sec:experiments} and summarize our findings in Section \ref{sec:discussion}.

{
We would like to point out that convergence of two time-scale methods is well-understood \cite{wang2016accelerating,wang2017stochastic}; however, applying these methods as is requires decision variables to be \emph{vectors}, not \emph{functions}. This parameterization, however, causes an approximation error which is difficult to characterize (see \cite{JMLR:v17:13-016}). In contrast, the RKHS parameterization, operating with the combination of projections and SQG, attains solutions that are close to the minimizer of the true Bellman evaluation error, where closeness is controlled by regularization introduced in the next section. This is due to the fact that RKHS possesses universal approximation under judicious choice of the kernel \cite{micchelli2006universal}, thus circumventing approximation error.}

{Recently, in companion work, an optimization-based variant of $Q$ learning in RKHS is developed \cite{tolstaya2018nonparametric}; however, a number of essential points distinguish that thread from methods developed here. Specifically, the optimization-based reformulation of Bellman's evaluation equation yields a convex program for which i.i.d. assumptions are close-to-valid. By contrast, Bellman's optimality equation yields a \emph{non-convex reformulation}. While the convergence of $Q$ learning requires i.i.d. assumptions, typically in practice these are violated. These statistical dependencies make policy learning a challenging domain to study Bayesian exploration, whereas policy evaluation is suitable \cite{precup2000eligibility}. Additionally, policy evaluation is just one component of reinforcement learning algorithms based upon policy search such as policy gradient method \cite{sutton2000policy} or actor-critic \cite{konda2003onactor,kumar2019sample}, whereas $Q$ learning is a standalone procedure. Other recent work focuses on policy search, a form of stochastic gradient with respect to a parameterized family of policies \cite{zhang2019global,8692920}, which is categorically different from fixed point iterations derived from Bellman's equations \cite{Bellman:1957}.}

\section{Policy Evaluation}\label{subsec:opt}
We reformulate the fixed point problem \eqref{eq:bellman_operator_v} defined by Bellman's equation so that it may be identified with a nested stochastic program. {Since the resulting domain of this problem is intractable, we hypothesize that the Bellman fixed point belongs to a RKHS. Then, to apply the Representer Theorem, we require the introduction of regularization.}

We proceed with reformulating \eqref{eq:bellman_operator_v}: subtract the value function $V^{\pi}(\bbx)$ that satisfies the fixed point relation from both sides, and then pull it inside the expectation:
\begin{equation}\label{eq:bellman_operator_v_reformulate1}
0 = \mathbb{E}_{\bby}[ r(\bbx,\pi(\bbx),\bby) + \gamma V^{\pi}(\bby) -  V^{\pi} (\bbx)    \given \bbx, \pi(\bbx)]
\end{equation} 
for all $\bbx \in \ccalX$. Value functions satisfying \eqref{eq:bellman_operator_v_reformulate1} are equivalent to those which satisfy the quadratic expression
%
$0=\frac{1}{2}( \mathbb{E}_{\bby}[ r(\bbx,\pi(\bbx),\bby) + \gamma V^{\pi}(\bby) -  V^{\pi}(\bbx)    \given \bbx, \pi(\bbx)])^2 \; ,$
which is null for all $\bbx \in \ccalX$, the starting point of the trajectory defining the value function \eqref{eq:value_function}. {To solve this expression for every $\bbx$, rather than solving it for a fixed $\bbx$ separately, we may integrate it out, which we do together with integrating over policy $\pi(\bbx)$ to yield the compositional stochastic program:}
\begin{align}
\label{eq:main_prob_0}
&V^{\pi} = \argmin_{V \in \ccalB(\ccalX)} J(V)  \\
&:= \!\argmin_{V \in \ccalB(\ccalX)} \mathbb{E}_{\bbx, \pi(\!\bbx\!)}\! \big\{\! \frac{1}{2}\!( \mathbb{E}_{\bby}\![ r(\bbx,\!\pi(\bbx\!)\!,\!\bby\!) \!\!+\! \gamma V\!(\bby\!) \!\!-\!\!  V \!(\!\bbx\!) \! \! \given\! \bbx,\! \pi(\!\bbx\!)])^2 \!\big\} \nonumber
\end{align}
whose solutions coincide exactly with the fixed points of \eqref{eq:bellman_operator_v}. The equivalence of  \eqref{eq:bellman_operator_v_reformulate1} and \eqref{eq:main_prob_0} is not in general true, but only true when the probability distribution $\mu$ over $\bbx$ is ergodic. That is, for fixed policy $\pi$, $\mu$ is non-vanishing over the entire state space $\ccalX$: for each $\bbx \in \ccalX$, $\mu(\bbx) > 0$. Henceforth, we require $\mu$, the prior distribution over states $\bbx\in\ccalX$, to be ergodic. {See \cite{bhatnagar2009convergent,sutton2009fast} for a discussions of transforming Bellman equations into objective functions, and the necessity of ergodicity.}

 \eqref{eq:main_prob_0} defines a functional optimization problem which is intractable when we search over all bounded continuous functions $\ccalB(\ccalX)${. However, when} we restrict $\ccalB(\ccalX)$ to a Hilbert space $\ccalH$ equipped with a unique \emph{reproducing kernel}, i.e., an inner product-like map $\kappa: \ccalX \times \ccalX \rightarrow \reals$ {such that for $f\in\ccalH$,}
\begin{align} \label{eq:rkhs_properties}
(i) \  \langle f , \kappa(\bbx, \cdot) \rangle _{\ccalH} \! = \! f(\bbx) \ 
(ii) \ \ccalH = \closure{\text{span}\{ \kappa(\bbx , \cdot) \}} \; ,
\end{align}
 %
for all $\bbx \in \ccalX$. We may apply the Representer Theorem to transform \eqref{eq:main_prob_0} into a semi-parametric one \cite{kimeldorf1971some,norkin2009stochastic}. In a RKHS, the optimal function $f\in\ccalH$ of \eqref{eq:main_prob_0} then takes the form \vspace{-1.5mm}
\begin{equation}\label{eq:kernel_expansion}
f(\bbx) = \sum_{n=1}^N w_n \kappa(\bbx_n, \bbx)\; {,}
\end{equation}
{where $\bbx_n$ is} a realization of {the} random variable $\bbx$. Thus, $f\in\ccalH$ is a kernel expansion \emph{only} at training samples. We define the upper summand index $N$  in \eqref{eq:kernel_expansion} in the kernel expansion of $f \in \ccalH$ {as the} model order, {which here} coincides with the training sample size.  
 Common kernel choices are {polynomials and radial basis  (Gaussian) functions}, i.e., $\kappa(\bbx,\bbx') = \left(\bbx^T\bbx'+b\right)^c $ and $\kappa(\bbx,\bbx') = \exp\{ -{\lVert \bbx - \bbx' \rVert_2^2}/{2c^2} \}$, respectively.
In \eqref{eq:rkhs_properties}, property (i) is called the reproducing property, which follows from Riesz Representation Theorem \cite{wheeden1977measure}. Replacing $f$ by $\kappa(\bbx' , \cdot) $  in \eqref{eq:rkhs_properties} (i) yields $ \langle \kappa(\bbx', \cdot) , \kappa(\bbx, \cdot) \rangle_{\ccalH} = \kappa(\bbx, \bbx')$, the origin of the term ``reproducing kernel."  Moreover, property \eqref{eq:rkhs_properties} (ii) states that functions $f\in\ccalH$ admit a basis expansion in terms of kernels  \eqref{eq:kernel_expansion}. 
Such spaces are {called reproducing kernel Hilbert spaces (RKHSs)}. {When the kernel is universal} \cite{micchelli2006universal}, e.g., a Gaussian, a continuous function over a compact set may be approximated uniformly by one in a RKHS.

Subsequently, we seek to solve {\eqref{eq:main_prob_0}} with $V \in \ccalH$, and independent and identically distributed samples $(\bbx_t, \pi(\bbx_t), \bby_t)$ from the triple $(\bbx, \pi(\bbx), \bby)$ are sequentially available, yielding
\begin{align}\label{eq:main_prob}
&V^{*}=\argmin_{V \in \ccalH} \mathbb{E}_{\bbx, \pi(\bbx)}\Big\{\frac{1}{2}( \mathbb{E}_{\bby}[ r(\bbx,\pi(\bbx),\bby)  \\
& \qquad\qquad\qquad\quad  + V(\bby) -  V (\bbx)    \given \bbx, \pi(\bbx)])^2\Big\} + \frac{\lambda}{2} \| V\|_{\ccalH}^2\nonumber
\end{align} 
Hereafter, define $L(V) :=  \mathbb{E}_{\bbx, \pi(\bbx)}\{\frac{1}{2}( \mathbb{E}_{\bby}[ r(\bbx,\pi(\bbx),\bby) + \gamma V(\bby) -  V (\bbx)    \given \bbx, \pi(\bbx)])^2\} $ and $J(V) = L(V) + (\lambda/2)\|V\|_{\ccalH}^2$. The regularization term $(\lambda/2)\|V\|_{\ccalH}^2$ in \eqref{eq:main_prob} is needed to apply the Representer Theorem \eqref{eq:kernel_expansion} \cite{kimeldorf1971some}. Thus, policy evaluation in infinite MDPs \eqref{eq:main_prob} is both a specialization of compositional stochastic programming \cite{wang2017stochastic} to an objective defined by dynamic programming, and a generalization to the case where the decision variable is not vector-valued but is instead a function.

\vspace{-1.5mm}
\section{Functional Stochastic Quasi-Gradient}
\label{sec:algorithm}
%
%
To apply functional SQG to \eqref{eq:main_prob}, {we} differentiate the compositional objective $L(V)$, which is of the form $L= g \circ h$, with $g(u)=\mathbb{E}_{\bbx, \pi(\bbx)} [(1/2) u^2]$ and $h(V) = \mathbb{E}_{\bby} [ r(\bbx, \pi(\bbx), \bby) + \gamma V(\bby) - V(\bbx) \given \bbx, \pi(\bbx)]$, and then consider its stochastic estimate. Consider the Frech$\acute{\text{e}}$t derivative of $L(V)$:
\begin{align}\label{eq:grad}
\nabla_V \! L(V)
%
%
&=\mathbb{E}_{\bbx, \pi(\bbx)}\big\{ \mathbb{E}_{\bby} [ \gamma\kappa(\bby, \cdot) - \kappa(\bbx,\cdot)\given \bbx, \pi(\bbx)] \\
&\qquad\quad \times\mathbb{E}_{\bby}\![ r(\bbx,\!\pi(\bbx),\!\bby\!) \!\!+\! \gamma V\!(\bby\!) \! - \! \!V \!(\bbx) \!   \given\! \bbx, \!\pi(\bbx)]  \big\} \nonumber
\end{align}
Here we pull the differential operator inside the expectation and use both the chain rule and reproducing property of the kernel \eqref{eq:rkhs_properties}(i). 
For future reference, we define the expression $\mathbb{E}_{\bby} [ r(\bbx, \pi(\bbx), \bby) + \gamma V(\bby) - V(\bbx) \given \bbx, \pi(\bbx)] = \bar{\delta}$ as the average temporal difference \cite{sutton1988learning}. To perform stochastic descent in function space $\ccalH$, we need a stochastic approximate of \eqref{eq:grad} evaluated at a state-action-state triple $(\bbx, \pi(\bbx), \bby)$, which together with the regularizer yields
\begin{align}\label{eq:stoch_grad}
{\nabla}_V &J(V , \delta ; \bbx, \pi(\bbx), \bby)  \\
&= \!\left[ \gamma\kappa(\bby, \!\cdot) \!-\! \kappa(\bbx,\!\cdot)\right] \!
\left[ r(\bbx,\pi(\bbx),\bby)\! +\! \gamma V\!(\bby\!) \!- \!  V \!(\bbx\!)\!   \right] \!+\! \lambda V\nonumber
\end{align}
where $\delta := r(\bbx, \pi(\bbx), \bby) + \gamma V(\bby) - V(\bbx)$ is defined as the (instantaneous) temporal difference.
Observe that we cannot obtain samples of ${\nabla}_V J(V, \delta ; \bbx, \pi(\bbx), \bby)$ with a single query to a simulation oracle: stochastic gradient method would estimate one of the expected gradients by its instantaneous approximation, but would still leave a second expected value that depends on infinitely many realizations of either prior distribution and policy $(\bbx, \pi(\bbx))$ or MDP transition dynamics $\bby$, a problem first identified in \cite{sutton2009convergent} for finite MDPs where it is called the \emph{double sampling problem}. {Double sampling aside, an equally significant  challenge associated with using \eqref{eq:stoch_grad} as a candidate descent direction is that classically we would compute its expectation conditional on the algorithm history, but due to the dependence of the factors, this does not yield \eqref{eq:grad}. In particular, the stochastic gradient is \emph{biased} with respect to \eqref{eq:grad} due to the inner conditional expectation in \eqref{eq:grad}.}

 To mitigate these issues, we {require} a method {that constructs} a \emph{coupled} stochastic descent procedure by considering noisy estimates of both factors in the product-of-expectations expression in \eqref{eq:grad}. 
The first factor $\left[ \gamma\kappa(\bby, \cdot) - \kappa(\bbx,\cdot)\right] $ in \eqref{eq:stoch_grad} is a difference of kernel maps, so {estimating its expectation is parameterized by infinitely many samples of $\bbx$ and $\bby$ \cite{Kivinen2004,smola2007hilbert}.} Instead, we propose a sequence based on samples of the second scalar factor to estimate its expected value. Specifically, from samples of $\delta
$, consider a recursion $z_t$ that estimates $\bar{\delta}$ as
\begin{align}\label{eq:td_average}
&\delta_t = r(\bbx_t, \pi(\bbx_t), \bby_t) + \gamma V_t(\bby_t) - V_t(\bbx_t)  \nonumber\\
&z_{t+1} = (1-\beta_t) z_t + \beta_t \delta_t
\end{align}
where we define $\delta_t$ \cite{sutton1988learning} as the temporal difference at time $t$ in \eqref{eq:td_average}.
Thus, \eqref{eq:td_average} averages the TD sequence $\delta_t$: $z_t$ estimates $\bar{\delta}_t$, and $\beta_t\in (0,1)$ is a learning rate.

To define a stochastic descent step, we replace the first factor inside the outer expectation in \eqref{eq:grad} with its instantaneous approximate, {i.e., $ \left[ \gamma\kappa(\bby_t, \cdot) - \kappa(\bbx_t,\cdot)\right]$,} at sample $(\bbx_t, \pi(\bbx_t), \bby_t)$, which yields the stochastic quasi-gradient step 
\begin{equation}\label{eq:quasi_fsgd}
\hat{V}_{t+1} =(1-\alpha_t\lambda)\hat{V}_t - \alpha_t (\gamma \kappa(\bby_t, \cdot) - \kappa(\bbx_t,\cdot))z_{t+1}  \; .
\end{equation}
where the coefficient $(1-\alpha_t \lambda)$ comes from the regularizer, and $\alpha_t$ is a positive scalar learning rate. This update is a stochastic quasi-gradient step because the true stochastic gradient of $J(V)$ is $(\gamma \kappa(\bby_t, \cdot) - \kappa(\bbx_t,\cdot))\delta_{t}$, but this estimator is unavailable with a single trajectory of the MDP since the factors in this product are dependent. By replacing $\delta_t$ by auxiliary variable $z_{t+1}$ this issue may be circumvented in the construction of coupled supermartingales (Section \ref{sec:convergence}).

{\bf Kernel Parameterization} Suppose $V_0=0\in \ccalH$. Then, \eqref{eq:quasi_fsgd} at time $t$, making use of the Representer Theorem \eqref{eq:kernel_expansion}, implies the function $\tilde{V}_t$ is a kernel expansion of past states $(\bbx_t, \bby_t)$ as\vspace{-1mm}
\begin{align}\label{eq:kernel_expansion_t}
\hat{V}_t(\bbx) 
= \sum_{n=1}^{2(t-1)} w_n \kappa(\bbv_n, \bbx)
= \bbw_t^T\boldsymbol{\kappa}_{\bbX_t}(\bbx) \; .
\end{align}
On the right-hand side of \eqref{eq:kernel_expansion_t} we introduce the notation $\bbv_n = \bbx_{n/2} \text{ for } n \text{ even}$ and $\bbv_n = \bby_{n/2 + 1} \text{ for } n \text{ odd}$, and:
%
%
$\bbw_t = [w_1, \cdots, w_{2(t-1)}] \in \reals^{2(t-1)} \; ,$ 
$\bbX_t = [\bbx_1, \bby_1, \ldots, \bbx_{t-1}, \bby_{t-1}]\in \reals^{p\times 2(t-1)} \; ,$
%
%
and $\boldsymbol{\kappa}_{\bbX_t}(\cdot) = [\kappa(\bbx_1,\cdot),\ \kappa(\bby_1,\cdot), \ \ldots\ ,\kappa(\bbx_{t-1},\cdot), \kappa(\bby_{t-1},\cdot)]^T \; .$
 The kernel expansion in \eqref{eq:kernel_expansion_t}, together with the functional update \eqref{eq:quasi_fsgd}, yields the fact that functional SQG in $\ccalH$ amounts to the following updates on the {data matrix $\bbX$, henceforth referred to as a kernel dictionary,} and coefficient vector $\bbw$:
\begin{align}\label{eq:param_update} 
&\bbX_{t+1} = [\bbX_t \; , \; \bbx_t \; , \; \bby_t] , \nonumber\\
&\bbw_{t+1} = [ (1 - \alpha_t \lambda) \bbw_t \; , \; \alpha_t z_{t+1} \; , \; -\alpha_t \gamma z_{t+1}]  \; ,
\end{align}
Observe that this update causes $\bbX_{t+1}$ to have two more columns than $\bbX_t$. We define the \emph{model order} as number of data points $M_t$ in the dictionary at time $t$, which for functional stochastic quasi-gradient descent is $M_t = 2(t-1)$. Asymptotically, then, the complexity of storing $\hat{V}_t(\bbx)$ is infinite. 

\subsection{Sparse Stochastic Subspace Projections}\label{subsec:polk_projection}
%

%
Since the update \eqref{eq:quasi_fsgd} has complexity $\ccalO(t)$ due to the RKHS parameterization \cite{Kivinen2004,koppel2019parsimonious}, it is impractical {in} settings with streaming data or arbitrarily large training sets. We address this issue by replacing the stochastic descent step \eqref{eq:quasi_fsgd} with an orthogonally projected variant \cite{koppel2019parsimonious}, where the projection is onto a low-dimensional functional subspace $\ccalH_{\bbD_{t+1}}$ of $\ccalH$, i.e.,
\begin{align}\label{eq:quasi_projected_fsgd}
\!\!\!\!\!\!V_{t+1}\! \!=\!\ccalP_{\ccalH_{\bbD_{t+1}}}\![(1 \!-\! \alpha_t \lambda\!) V_t \!  -\! \alpha_t (\!\gamma \kappa(\bby_t, \!\cdot) \!-\! \kappa(\!\bbx_t,\cdot\!))z_{t+1}\! ]  ,
\end{align}
where $\alpha_t$ again is a scalar step-size, and $\ccalH_{\bbD_{t+1}}=\text{span}\{\kappa(\bbd_n, \cdot) \}_{n=1}^{M_t}$ for some collection of sample instances $\{\bbd_n\}\subset \{\bbx_u\}_{u\leq t}$. Note that the un-projected function SQG method \eqref{eq:quasi_fsgd} may be interpreted as conducting a sequence of orthogonal projections, which motivates the design of \eqref{eq:quasi_projected_fsgd}. Specifically, rewrite \eqref{eq:quasi_fsgd} as the quadratic minimization
\begin{align}\label{eq:proximal_hilbert_dictionary}
\hat{V}_{t+1} &\!\!=\! \argmin_{V \in \ccalH} \!
\Big\lVert  V \!\!\!-\! \!\Big(\!\!(\!1\!-\!\alpha_t \lambda\!)\! \hat{V}_t \!
\!-\! \alpha_t \!  (\!\gamma \kappa(\!\bby_t,\! \cdot\!) \!
 - \!\kappa(\!\bbx_t,\!\cdot\!)\!)z_{t+1}\!\! \Big)\!\Big\rVert_{\ccalH}^2 \nonumber \\
&\!\!=\!\! \argmin_{V \in \ccalH_{\bbX_{t+1}}}\!\! \!
\Big\lVert  V \!\!\!-\! \!\Big(\!\!(\!1\!-\!\alpha_t \lambda\!)\! \hat{V}_t \!
\!-\! \alpha_t \!  (\!\gamma \kappa(\!\bby_t,\! \cdot\!) \!
 - \!\kappa(\!\bbx_t,\!\cdot\!)\!)z_{t+1}\!\! \Big)\!\Big\rVert_{\ccalH}^2, 
\end{align}
where the first equality in \eqref{eq:proximal_hilbert_dictionary} comes from ignoring constant terms which vanish upon differentiation with respect to $V$, and the second comes from observing that $V_{t+1}$ can be represented using only the points $\bbX_{t+1}$, using \eqref{eq:param_update}.  Notice \eqref{eq:proximal_hilbert_dictionary} expresses $V_{t+1}$ as the projection $(1-\alpha_t \lambda) V_t 
- \alpha_t   (\gamma \kappa(\bby_t, \cdot) - \kappa(\bbx_t,\cdot))z_{t+1}$ onto the subspace defined by dictionary $\bbX_{t+1}$.

Rather than {selecting} dictionary $\bbD=\bbX_{t+1}$, we propose instead to select a different dictionary, $\bbD=\bbD_{t+1}$, which is extracted from the data points observed thus far, at each iteration.  {The process by which we select $\bbD_{t+1}$ is delayed for now, but is of dimension $p \times {M}_{t+1}$. We design a scheme such that  ${M}_{t+1}$ is \emph{independent} of $t$, and instead determined by fundamental topological properties of state space $\ccalX$, i.e., a generalization of the Nyquist rate \cite{zhou2002covering}.}  As a result, the sequence ${V}_t$ differs from the functional stochastic quasi-gradient method $\hat{V}_t$ presented at the outset of this section. 

Specifically, suppose the function ${V}_{t+1}$ is parameterized dictionary $\bbD_{t+1}$ and weight vector $\bbw_{t+1}$. We denote columns of $\bbD_{t+1}$ as $\bbd_n$ for $n=1,\dots,{M}_{t+1}$, where the time index is dropped for notational clarity but may be inferred from the context.
Setting aside how $\bbD_{t+1}$ is chosen for now, we replace the update \eqref{eq:proximal_hilbert_dictionary} in which the dictionary grows at each iteration by the  functional stochastic quasi-gradient sequence projected onto the subspace $\ccalH_{\bbD_{t+1}}=\text{span}\{ \kappa(\bbd_n, \cdot) \}_{n=1}^{M_{t+1}}$ as\vspace{-2mm}
\begin{align}\label{eq:projection_hat}
{V}_{t+1}&\!=\!\! \argmin_{V \in \ccalH_{\bbD_{t+1}}}\!\! \!
\Big\lVert  V \!\!\!-\! \!\Big(\!\!(\!1\!-\!\alpha_t \lambda\!)\! \hat{V}_t \!
\!-\! \alpha_t \!  (\!\gamma \kappa(\!\bby_t,\! \cdot\!) \!
 - \!\kappa(\!\bbx_t,\!\cdot\!)\!)z_{t+1}\!\! \Big)\!\Big\rVert_{\ccalH}^2, \nonumber \\
&:=\!\ccalP_{\ccalH_{\bbD_{t+1}}} \!\!\Big[\! (1\! - \!\alpha_t \lambda) V_t\! -\! \alpha_t (\gamma \kappa(\bby_t, \cdot) \! -\! \kappa(\bbx_t,\cdot))z_{t+1} \Big] .
\end{align}
where we define the projection operator $\ccalP$ onto subspace $\ccalH_{\bbD_{t+1}}\subset \ccalH$ by the update \eqref{eq:projection_hat}. This orthogonal projection is the modification of the functional SQG iterate [cf. \eqref{eq:quasi_fsgd}] defined at the beginning of this subsection \eqref{eq:quasi_projected_fsgd}. Next we discuss how this update amounts to modifications of the parametric updates \eqref{eq:param_update} defined by functional SQG.
These subspace projections may be computed efficiently by exploiting the kernel parameterization described in Appendix \ref{apx_parameterization} operating together with destructive matching pursuit \cite{Vincent2002}.


%
\begin{algorithm}[t]
\caption{PKGTD: Parsimonious Kernel Gradient Temporal Difference}
\begin{algorithmic}
\label{alg:pkgtd}
\REQUIRE $\{\bbx_t,\pi(\bbx_t), \bby_t,\alpha_t,\beta_t,\eps_t \}_{t=0,1,2,...}$
\STATE \textbf{initialize} ${V}_0(\cdot) = 0, \bbD_0 = [], \bbw_0 = [], z_0=0$, i.e. initial dict., coeffs., and aux. variable null
\FOR{$t=0,1,2,\ldots$}
	\STATE Obtain trajectory realization $(\bbx_t,\pi(\bbx_t), \bby_t)$
	\STATE Compute the TD and auxiliary sequence $z_{t+1}$  [cf. \eqref{eq:td_average}]:\vspace{-2mm}
$$\quad \delta_t = r(\bbx_t, \pi(\bbx_t), \bby_t) + \gamma V_t(\bby_t) - V_t(\bbx_t),$$
$$z_{t+1} = (1-\beta_t) z_t + \beta_t \delta_t$$\vspace{-4.5mm}
	\STATE Compute unconstrained functional SQG step [cf. \eqref{eq:quasi_fsgd}]\vspace{-2mm}
	$$\tilde{V}_{t+1}(\cdot) =(1-\alpha_t \lambda)\tilde{V}_t(\cdot) - \alpha_t (\gamma \kappa(\bby_t, \cdot) - \kappa(\bbx_t,\cdot))z_{t+1}$$\vspace{-4.5mm}
	\STATE Revise dict. $\tbD_{t+1} = [\bbD_t,\;\! \bbx_t \; , \bby_t ]$, weights $\tbw_{t+1} \leftarrow [(1-\alpha_t\lambda)\bbw_t,\;\;\alpha_t z_{t+1} ,-\alpha_t \gamma z_{t+1}]$
	\STATE Compress dictionary via Alg. \ref{alg:komp}, obtain coeffs. via \eqref{eq:proximal_hilbert_representer} 
	$$({V}_{t+1},\bbD_{t+1},\bbw_{t+1}) = \textbf{KOMP}(\tilde{V}_{t+1},\tbD_{t+1},\tbw_{t+1},\epsilon_t)$$\vspace{-3mm}
\ENDFOR
\end{algorithmic}
\end{algorithm}


We summarize the overall method, Parsimonious Kernel Gradient Temporal Difference (PKGTD) in Algorithm \ref{alg:pkgtd}: we execute the stochastic projection of the functional SQG iterates onto sparse subspaces $\ccalH_{\bbD_{t+1}}$ stated in \eqref{eq:projection_hat}. With initial function null $V_0=0$ (empty dictionary $\bbD_0=[]$ and coefficients $\bbw_0=[]$),
at each step, given an i.i.d. sample $(\bbx_t, \pi(\bbx_t),\bby_t)$ and step-sizes $\alpha_t, \beta_t$, we compute the \emph{unconstrained} functional SQG iterate $\tilde{V}_{t+1}(\cdot) = (1-\alpha_t \lambda)\tilde{V}_t(\cdot) - \alpha_t (\gamma \kappa(\bby_t, \cdot) - \kappa(\bbx_t,\cdot))z_{t+1}$ parameterized by $\tbD_{t+1}$ and $\tbw_{t+1}$ as stated in \eqref{eq:param_tilde}, which are fed into KOMP (Algorithm \ref{alg:komp} in Appendix \ref{apx_parameterization}) with budget $\eps_t$, i.e., $(V_{t+1}, \bbD_{t+1}, \bbw_{t+1})= \text{KOMP}(\tilde{V}_{t+1},\tilde{\bbD}_{t+1}, \tilde{\bbw}_{t+1},\eps_t)$.

{
\begin{remark}\label{remark_twotimescale}\normalfont
While two time-scale stochastic approximation originally appeared in the 1980s \cite{Korostelev,ermoliev1983stochastic} for compositional stochastic programming with asymptotic stability established, their role in RL, namely, to form the foundation of actor-critic algorithms \cite{borkar1997actor,konda2003onactor} was more recent. 

A separate but related line of research in RL identifies their use not in actor-critic (which is at its core a policy search method), but instead in order to solve Bellman equations (approximate dynamic programming), beginning with \cite{sutton2009convergent}. In \cite{sutton2009convergent} the authors hypothesize that one possible reason for instability of temporal difference learning under function approximation (as detailed in \cite{tsitsiklis1997analysis}), is that these are \emph{not gradient algorithms} but instead stochastic fixed point iterations. Thus, their stability interacts in a more intricate manner with the function parameterization. By reformulating Bellman equations as optimization problems, these problems are identified as possessing compositional structure, and thus are amenable to two time-scale algorithms, yielding \emph{gradient temporal difference} learning (GTD). The derivation of PKGTD is structurally aligned with GTD in that its derivation generalizes the derivation of GTD, rather than actor-critic, although two time-scale methods are at the core of both approaches. 
\end{remark}}




%
\section{Convergence Analysis}\label{sec:convergence}
{We now analyze the stability and memory requirements of Algorithm \ref{alg:pkgtd} developed in Section \ref{sec:algorithm}. In stochastic fixed-point methods such as TD learning{,} the interplay between the Bellman operator contraction \cite{bertsekas1978stochastic} and expectations prevents the construction of supermartingales underlying stochastic descent stability \cite{Robbins1951}. {Attempts to overcome this challenge based on stochastic backward-differences require the state space to be completely explored in the limit \emph{per step} (intractable when $|\ccalX| = \infty$) \cite{tsitsiklis1994asynchronous}, or stipulate that data dependent matrices be non-singular \cite{sutton2009convergent}, respectively.} Thus these methods must be analyzed using ideas from dynamical systems \cite{borkar2000ode}. 
 In contrast, we establish that Algorithm \ref{alg:pkgtd} belongs to the family of descent algorithms, and hence its behavior can be connected to that of supermartingales \cite{wang2013incremental} -- to the best of our knowledge, this is the first time supermartingales have been used in analyzing stochastic methods for MDPs. This is also true of GTD \cite{sutton2009convergent}, although it is analyzed using ODEs \cite{borkar2000ode}. 
}

%

Under the assumptions stated next, it is possible to derive the fact that the auxiliary variable $z_t$ and value function estimate $V_t$ satisfy supermartingale-type relationships, but their behavior is intrinsically coupled. We generalize recently developed coupled supermartingale tools in \cite{wang2013incremental}, i.e., Lemma \ref{lemma:coupled_supermartingale} in Appendix \ref{apx_assumptions}, to RKHSs in order to establish the following almost sure convergence result when the step-sizes and compression budget are diminishing. More specifically, we require the following technical definitions and conditions.

\subsection{Technical Assumptions and Definitions}\label{apx_assumptions1}

In particular, for further reference, we define the functional stochastic quasi-gradient of the regularized objective as
\begin{align}\label{eq:quasi_sg}
&\hat{\nabla}_V J(V_t, z_{t+1}; \bbx_t, \pi(\bbx_t), \bby_t )=\nonumber\\
&\qquad\qquad\qquad (\gamma \kappa(\bby_t, \cdot) - \kappa(\bbx_t,\cdot))z_{t+1}  + \lambda V_t \; ,
\end{align}
and its sparse-subspace projected variant as
\begin{align}\label{eq:projected_quasi_sg}
&\tilde{\nabla}_V J(V_t, z_{t+1}; \bbx_t, \pi(\bbx_t), \bby_t ) \\
&= \!\frac{1}{\alpha_t} \Big( \!V_{t} \!-\!  \ccalP_{ \ccalH_{\bbD_{t+1}}} \Big[ V_t \!-\! \alpha_t \hat{\nabla}_V J(V_t, z_{t+1}; \bbx_t, \pi(\bbx_t), \bby_t )\!\Big]\!\Big)  \; ,\nonumber
\end{align}
Note that the update \eqref{eq:quasi_projected_fsgd}, using \eqref{eq:projected_quasi_sg}, may be rewritten as a stochastic projected quasi-gradient step rather than a stochastic quasi-gradient step followed by set projection, i.e.,
\begin{equation}\label{eq:quasi_projected_fsgd_projected}
V_{t+1} =V_t - \alpha_t \tilde{\nabla}_V J(V_t, z_{t+1}; \bbx_t, \pi(\bbx_t), \bby_t )  \; ,
\end{equation}

{Further, define the time-dependent sigma algebra, i.e., filtration, as  $\ccalF_t \supset (\{V_s, z_s,\bbx_s, \pi(\bbx_s), \bby_s\}_{s=0}^{t-1})$.} Now we are ready to state the technical conditions required for convergence. All statements involving conditional expectations are imposed with probability $1$, unless otherwise stated.

\begin{assumption}\label{as:first}
The state space $\ccalX\subset\reals^p$ and action space $\ccalA\subset \reals^q$ are compact, and the reproducing kernel map may be bounded as
\begin{equation}\label{eq:bounded_kernel}
\sup_{\bbx\in\ccalX} \sqrt{\kappa(\bbx, \bbx )} = X < \infty
\end{equation}
%
\end{assumption}

\begin{assumption}\label{as:mean_variance_td} 
The temporal difference $\delta$ and auxiliary sequence $z$ [cf. \eqref{eq:td_average}] satisfy the zero-mean and finite conditional variance conditions, respectively,
\begin{align}\label{eq:td_assumption}
&\mathbb{E}\left[\delta \given \bbx, \pi(\bbx) \right] =\bar{\delta} \; ,
 \qquad\quad \mathbb{E}\left[(\delta - \bar{\delta})^2 \given \ccalF_t \right] \leq \sigma_{\delta}^2 \; , \nonumber\\
& \mathbb{E}\left[z^2 \given \bbx, \pi(\bbx) \right] \leq G_{\delta}^2\; .
\end{align}
where  $\sigma_{\delta}$ and $G_{\delta}$ are positive scalars. 
\end{assumption}

\begin{assumption}\label{as:mean_variance_frechet} 
The stochastic quasi-gradient, when evaluated at $\bar{\delta}$, is an unbiased estimate for the true gradient $\nabla_V J(V)$. Moreover, the difference of reproducing kernels expression (the first factor in \eqref{eq:stoch_grad}) has finite conditional variance:
\begin{align}\label{eq:value_kernel_assumption}
&\mathbb{E}\left[ (\gamma\kappa(\bby, \cdot) - \kappa(\bbx,\cdot)) \bar{\delta} \right] = \nabla_V J(V) \; , \nonumber\\
& \mathbb{E}\left[ \| \gamma\kappa(\bby_t, \cdot) - \kappa(\bbx_t,\cdot)\|_{\ccalH}^2 \given \ccalF_t \right ] \leq G_{V}^2 \; . 
\end{align}
Additionally, the projected stochastic gradient of the objective [cf. \eqref{eq:projected_quasi_sg}] has finite second conditional moment as
\begin{equation}\label{eq:projected_quasi_sg_variance}
\mathbb{E}\left[\|\tilde{\nabla}_V J(V_t, z_{t+1}; \bbx_t, \pi(\bbx_t), \bby_t ) \|_{\ccalH}^2 \given \ccalF_t \right] \leq \sigma_V^2 \; ,
\end{equation}
and the conditional mean of the temporal difference $\bar{\delta}$ is Lipschitz continuous with respect to the value function $V$, i.e for any two distinct $\delta$ and $\tilde{\delta}$, we have
\begin{align}\label{eq:td_lipschitz}
|\bar{\delta} - \bar{\tilde{\delta}}  | \leq L_{V} \| V - \tilde{V} \|_{\ccalH}
\end{align}
where $V,\tilde{V} \in \ccalH$ are distinct RKHS elements, $L_V >0$ is a scalar, and $ \bar{\delta} = \mathbb{E}_{\bby} [ r(\bbx, \pi(\bbx), \bby) + \gamma V(\bby) - V(\bbx) \given \bbx, \pi(\bbx)] $.
\end{assumption}

Assumption \ref{as:first} regarding the compactness of the state and action spaces of the Markov Decision Process intrinsically hold for most application settings and limit the radius of the set from which the MDP trajectory is sampled. Similar boundedness conditions on the reproducing kernel map have been considered in supervised learning applications \cite{Kivinen2004}. The mean and variance properties of the temporal difference stated in Assumption \ref{as:mean_variance_td} to bound the error in the descent direction associated with stochastic approximations, and are necessary to establish stability of stochastic methods. Assumption \ref{as:mean_variance_frechet} is similar to Assumption \ref{as:mean_variance_td} but instead of establishing bounds on the stochastic approximation error of the temporal difference, limits stochastic error variance in the reproducing kernel Hilbert space. These are natural extensions of the conditions needed for convergence of stochastic compositional gradient methods with vector-valued decision variables \cite{wang2017stochastic}. {However, we note that \eqref{eq:td_lipschitz}, in the context of MDPs, restricts the class of reward functions to be those which may be smoothly interpolated in a RKHS.} {This condition holds, for instance, when the reward is a potential or navigation-like function \cite{rimon1992exact,ngicml1999}, which are well-known to interact favorably in the design of controllers from a dynamical systems perspective.}

{The stipulation that the KOMP projection explicitly thresholds the norm of the value functions allows us to write}
\begin{equation}\label{eq:bounded_iterates}
\| V_t \|_{\ccalH} \leq K \; , \qquad \|V^* \|_{\ccalH} \leq K \; , \quad \text{ for all } t
\end{equation}
where $K>0$ is some constant. The boundedness of $V^*$ follows from the fact that since $\ccalX$ is compact and $J(V)$ is a continuous convex function over a compact set, its minimizer is achieved over this compact set \cite{brezis2010functional}[Corrolary 3.23].

\begin{theorem}\label{theorem:convergence_wp1}
Consider the sequence $z_t$ [cf. \eqref{eq:td_average}] and $\{V_t\}$ [cf. \ref{eq:quasi_projected_fsgd}] as stated in Algorithm \ref{alg:pkgtd}. Assume the regularizer is positive $\lambda > 0$, Assumptions \ref{as:first} - \ref{as:mean_variance_frechet} hold, with the step-size conditions:
%
%
%
%
\begin{align}\label{eq:step_size_condition}
\!\!\!\sum_{t=1}^\infty\! \alpha_{t} \!=\!\infty \; ,\   \sum_{t=1}^\infty\! \beta_t \!= \!\infty , \  
\sum_{t=1}^\infty\! \alpha_{t}^2 \!+ \! \beta_t^2\!+\! \frac{\alpha_t^2}{\beta_t} \!\! <\!\! \infty\;, \ \eps_t \!= \!\alpha_t^2
\end{align}
Then $V_t \rightarrow V^*$ [cf. \eqref{eq:main_prob}] with probability $1$, and thus achieves the regularized Bellman fixed point \eqref{eq:bellman_operator_v_reformulate1} restricted to the RKHS.
\end{theorem}
%

The proof is given in Appendix \ref{apx_theorem_convergence_wp1}. Theorem \ref{theorem:convergence_wp1} states that the value functions generated by Algorithm \ref{alg:pkgtd} converge a.s. to the optimal $V^*$ defined by \eqref{eq:main_prob}. {With regularizer $\lambda$ made small but nonzero, using a universal kernel (e.g., a Gaussian), $V_t$ converges close to a function satisfying Bellman's equation in \emph{infinite MDPs} \eqref{eq:bellman_operator_v}. By decreasing the regularizer, limiting solutions close in on those which satisfy Bellman's  equation, though precise notions of closeness require continuity which is difficult to verify, given an arbitrary bounded reward. }This is the first guarantee w.p.1 for a true stochastic descent method with an infinitely and nonlinearly parameterized value function.  
Theorem  \ref{theorem:convergence_wp1} requires attenuating step-sizes such that the stochastic approximation error approaches null. In contrast, constant learning rates allow for maintaining algorithm adaptivity, motivating the following result.

One step-size sequence which satisfies the attenuation conditions \eqref{eq:step_size_condition} is 
%
%
$\alpha_t =\ccalO( t^{-(3/4 +\zeta/2)} )\; , \beta_t = \ccalO(t^{-(1 +\zeta)/2 })\; , \eps_t=\ccalO(\alpha_t^2) = \ccalO(t^{-(3/2 + \zeta)})$, where $\zeta>0$ is an arbitrarily small constant so that series $\sum_t \alpha_t$ and $\sum_t \beta_t$ diverge. Generally, satisfying \eqref{eq:step_size_condition}, requires: $\alpha_t = \ccalO(t^{-p_{\alpha}})$, $\beta_t = \ccalO(t^{-p_{\beta}})$ with $p_{\alpha}\in(3/4, 1)$ and $p_{\beta}\in (1/2, 2 p_{\alpha}- 1)$.

\begin{theorem}\label{theorem:constant_stepsize_convergence}
Suppose Algorithm \ref{alg:pkgtd} is run with constant positive learning rates $\alpha_t=\alpha$ and $\beta_t=\beta$ and constant compression budget $\eps_t=\eps$ with sufficiently large regularization, i.e.
\begin{align}\label{eq:constant_stepsize_parameters}
0<\beta < 1 \; , \alpha =\beta, \eps=C \alpha^2, \lambda=G_V^2 + \lambda_0 
\end{align}
where $C>0$ is a scalar, and $0<\lambda_0 < 1$. Then, under Assumptions \ref{as:first} - \ref{as:mean_variance_frechet}, the sub-optimality sequence $\|V_t - V^*\|_{\ccalH}^2$ converges in mean to a neighborhood:
\begin{align}\label{eq:constant_stepsize_convergence}
\! \! \liminf_{t\rightarrow \infty} \E{\|V_{t}\! -\! V^*\|_{\ccalH}^2}={\ccalO\left(\alpha  \right)} \; .
\end{align}
\end{theorem}
Theorem \ref{theorem:constant_stepsize_convergence} (proof in Appendix \ref{apx_theorem_constant_stepsize_convergence}) establishes that the value function estimates generated by Algorithm \ref{alg:pkgtd} converge in expectation to a neighborhood when constant step-sizes $\alpha$ and $\beta$ and sparsification budget $\epsilon$ in Algorithm \ref{alg:komp} are small constants. In particular, the bias $\eps$ induced by sparsification does not cause instability even when it is \emph{not going to null}. Moreover, this result only holds when the regularizer $\lambda$ is chosen large enough, which numerically induces a forgetting factor on past kernel dictionary weights \eqref{eq:param_tilde}. 
We may make the learning rates $\alpha$ and $\beta$ arbitrarily small, which yield a proportional decrease in the radius of convergence to a neighborhood of the Bellman fixed point \eqref{eq:bellman_operator_v}. 

{In general, Theorem \ref{theorem:constant_stepsize_convergence} \emph{does not} imply the sequence actually converges, but only that its $\liminf$ converges. To establish that the entire sequence converges, even in expectation, when used with constant step-sizes, one must recursively average the  functions and apply convexity of the objective function to establish the error bound decreases with the final iteration, as in Polyak-Ruppert averaging \cite{polyak1992acceleration}. Averaging, however, will be afflicted by the fact that different functions in the RKHS do not belong to the same subspace, and therefore their kernel dictionaries will need to be pooled, causing the model order to spike. Thus, averaging is a technique  of theoretical interest only in establishing limiting genuine behavior in RKHSs, and cannot be used unless parsimony is not a consideration. 
%
%
%
By contrast, under constant step-sizes selection, the value function estimates have moderate complexity in the worst case.

\begin{remark}\label{remark1} \normalfont
(Aggressive Constant Learning Rates)
In practice, one may obtain better performance by using larger constant step-sizes. To do so, the criterion \eqref{eq:constant_stepsize_parameters} may be relaxed: we require $0<\beta<1$ but $\alpha>0$ may be any positive scalar. Then, with regularizer chosen as $\lambda=G_V^2\frac{\alpha}{\beta} + \lambda_0 $ for $0<\lambda_0<1$, the radius of convergence is (see Appendix \ref{apx_theorem_constant_stepsize_convergence})
\begin{align}\label{eq:constant_stepsize_convergence2}
&{\liminf_{t\rightarrow \infty} \E{\|V_{t} - V^*\|_{\ccalH}^2} }\nonumber\\
&\qquad =\ccalO\Big(\alpha^2 + \beta^2 + \frac{\alpha^2}{\beta}\Big[1 + \alpha^2 + \frac{\alpha}{\beta} + \frac{\alpha^2}{\beta^2}\Big] \Big) \; .
\end{align}
The ratios $\alpha^2/\beta$ and $\alpha^2/\beta^2$ dominate \eqref{eq:constant_stepsize_convergence2} and must be made small  to obtain accurate solutions. 
\end{remark}
%
\begin{remark}\label{remark2} \normalfont
(Regularization Path)
{In Theorem \ref{theorem:convergence_wp1}, we establish convergence for any $\lambda>0$ when step-sizes attenuate. That regularizer $\lambda$ may not be null means that we do not extract the exact Bellman fixed point restricted to the RKHS, but only a function that is close. In related work the minimizer $V_{\lambda}^*$ continuously depends on $\lambda$. It's beyond the scope of this work to extend these results to this setting, but on the hypothesis that they generalize to \eqref{eq:main_prob}, we may claim that decreasing $\lambda$ causes $V_{\lambda}^*$ to be closer to fixed point stated in \eqref{eq:bellman_operator_v}. 

On the other hand, for Theorem \ref{theorem:constant_stepsize_convergence} with given regularizer $\lambda$, imposes explicit restrictions on the choice of constant algorithm step-sizes. That is, we require $\lambda=G_V^2\frac{\alpha}{\beta} + \lambda_0 $ for $0<\lambda_0<1$, where $\alpha > 0$ and $0<\beta<1$. It is possible to derive the fact that larger regularization means faster learning rates but to less accurate solutions, in either diminishing or constant step-size settings, but these facts, which depend on rate analyses, are left to future work.}
\end{remark}\vspace{-3mm}
%

As noted in Section \ref{sec:algorithm}, the complexity of functional stochastic quasi-gradient method in a RKHS is of order $\ccalO(2(t-1))$ which grows without bound. {To surmount this challenge, we propose subspace projections in Section \ref{subsec:polk_projection}.}  We {formalize here} that this projection indeed controls complexity when constant learning rates and compression budget are used. This result is a corollary, as it extends Theorem 3 in \cite{koppel2019parsimonious}. To obtain this result (proof in Appendix \ref{apx_corollary1}), we require the reward function to be bounded, as we state next.


\begin{assumption}\label{as:reward} 
The reward $r: \ccalX \times \ccalA \times \ccalX \rightarrow \reals$ is bounded:
\begin{equation}\label{eq:bounded_reward}
r(\bbx_t, \pi(\bbx_t), \bby_t) \leq R_{\max} \; \text{ for all } t, \bbx, \bba, \bby
\end{equation}
\end{assumption}

Assumption \ref{as:reward} holds whenever the reward function is continuous and the state and action spaces are compact, and thus holds for many popular RL problems. In this setting, the complexity of Algorithm \ref{alg:pkgtd} is finite, as is formalized next.

\begin{figure*}[t]
\centering
\vspace{-5mm}
  \includegraphics[width=10cm,height=1cm]{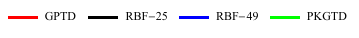}\\
  %
\hspace{1mm} \includegraphics[width=5.75cm,height=5cm]{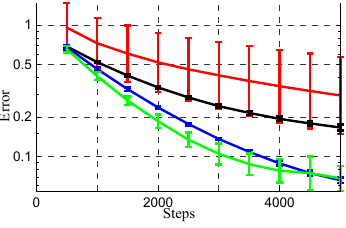}\hspace{2mm}
  %
    \includegraphics[width=5.75cm,height=5cm]{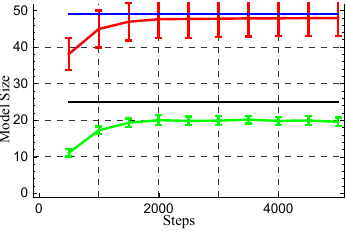}\hspace{2mm}
      %
  \includegraphics[width=5.8cm,height=5.1cm]{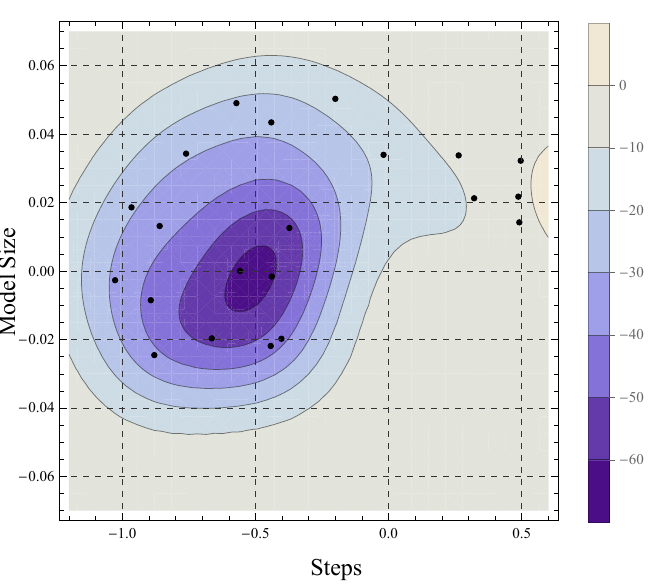}
  \caption{Experimental comparison of PKGTD to existing kernel methods for policy evaluation on the Mountain Car task. Test set  error (left), and the parameterization complexity (center) vs. iterations. PKGTD learns fastest and most stably with the least complexity (best viewed in color). We plot the contour of the learned value function (right): its minimal value is in the valley, and states near the goal are close to null. Bold black dots are kernel dictionary elements, or retained instances. {Sample means and standard deviations of performance metrics (left and center) are generated via $100$ individual training runs.} }
  \label{fig:experiments}\vspace{-2mm}
\end{figure*}

\begin{corollary}\label{corollary1}
Denote $V_t$ as the value function of Algorithm \ref{alg:pkgtd} with constant step-sizes $\alpha_t = \alpha$ and $\beta_t = \beta \in (0, 1)$ with compression budget $\eps_t =  \eps =C \alpha^2 $ and regularization $\lambda=(\alpha/\beta) G_V^2 + \lambda_0= \ccalO(\alpha\beta^{-1}+1) $  as in Remark \ref{remark1}. Let $M_t$ be its associated model order, i.e., the number of columns in its dictionary. Then there exists a finite upper bound $M^\infty$ such that, for all $t\geq0$, the model order is bounded $M_t\leq M^\infty$.
\end{corollary}
%

%
%
%
{Under specific selections \eqref{eq:constant_stepsize_parameters}, the {algorithm converges to a neighborhood} of the optimal {value function, whose radius depends} on the step-sizes, and may be made small by decreasing $\alpha$ at the cost of a decreasing learning rate. More importantly, {the use} of constant step-sizes and compression budget with large enough regularization yields a value function parameterized by a dictionary whose model order is always bounded (Corollary \ref{corollary1}). Thus, we may converge to an optimal neighborhood while ensuring the memory of the function parameterization is under control, and in the worst-case related to the covering number (metric entropy) of the state space.}


\vspace{-2mm}
\section{Experiments}
\label{sec:experiments}
%
%
%
%
%
%
%

Our experiments aim to compare PKGTD to other policy evaluation techniques in this domain.  Because it seeks memory-efficient solutions over an RKHS, we expect PKGTD to obtain accurate estimates of the value function using only a fraction of the memory required by the other methods.
We perform experiments on the classical Mountain Car domain \cite{sutton1998reinforcement}: an agent applies discrete actions $\ccalA = \left\{ \texttt{reverse}, \texttt{coast}, \texttt{forward} \right\}$ to a car that starts at the bottom of a valley and attempts to climb up to a goal at the top of one of the mountain sides. The state space is continuous, consisting of the car's scalar position and velocity, i.e., {$\ccalX \subset \reals^2$}.  The reward function $r(\bbx_t,\bba_t,\bby_t)$ is $-1$ unless $\bby_t$ is the goal state at the mountain top, in which case it is $0$ and the episode terminates.

Now we describe the configuration of the algorithms used for comparison. The Mountain Car environment has a two-dimensional state space (position and velocity) with bounds of $[-1.2, 0.6]$ in position, and $[-0.07, 0.07]$ in velocity. We chose not to normalize this state space to $[0,1]$ intervals, choosing instead to handle the scale difference by using non-isotropic kernels. The ratio of the kernel variances is equal to the ratio of the lengths of their corresponding bounds, so they would be isotropic kernels if we normalized the state space.

We used a fixed non-isotropic kernel bandwidth of $\sigma_1 = 0.2, \sigma_2 = 0.0156$ in all cases. By fixing the kernel bandwidth across all algorithms, we are basically enforcing that the learned functions all belong to the same Kernel Hilbert Space.

%

For PKGTD, the relevant parameters are the step size, $\alpha$, the rate of expectation update, $\beta$, the regularizer, $\lambda$, and the approximation error, $K$. For GPTD, the relevant parameters are the gaussian process noise standard deviation, $\sigma_0$, the linear independence test bound, $\nu$, and the regularizer, $\lambda$. For the RBF grids fit using GTD, the relevant parameters are the grid spacing in the position and velocity directions, $h_1$ and $h_2$, respectively, the step size, $\alpha$, and the rate of expectation update, $\beta$. Our values are summarized in Table 1. {While theoretically the selection of $\lambda$ is sensitive, experimentally we find it to have little impact, and thus fix it as a small value $\lambda = 10^{-6}$. 
\begin{tabular}{lcccccccc}
           &   \hspace{-2mm} $\alpha$ &\hspace{-6mm} $\beta$ & \hspace{-6mm} $\lambda$ & \hspace{-2mm} $K$   &  \hspace{-2mm} $\sigma_0$ &\hspace{-2mm} $\nu$ &\hspace{-2mm} $h_1$ & \hspace{-2mm}$h_2$  \\
  \cline{2-9}
  \hspace{-3mm} PKGTD    &  \hspace{-4mm} 8.0      & \hspace{-4mm}0.2     &  \hspace{-2mm}1e-6      &  \hspace{-2mm}0.02  &    \hspace{-2mm}        &       &       &        \\
\hspace{-3mm}  GPTD     &     \hspace{-4mm}      &     \hspace{-4mm}    &  \hspace{-2mm}1e-6      &     \hspace{-2mm}   &\hspace{-2mm} 0.01       & 0.2   &       &        \\
\hspace{-3mm}  RBF-25   & \hspace{-4mm} 10.0     &\hspace{-4mm} 0.25    &     \hspace{-2mm}       &    \hspace{-2mm}    &  \hspace{-2mm}          & \hspace{-2mm}      &\hspace{-2mm} 0.44  &\hspace{-2mm} 0.0343 \\
\hspace{-3mm}  RBF-49   & \hspace{-5.5mm} 1.5      & \hspace{-3mm}0.35    &       \hspace{-2mm}     &    \hspace{-2mm}    &    \hspace{-2mm}        &    \hspace{-2mm}   & \hspace{-1mm}0.26  & \hspace{-1mm}0.0203 \\
\end{tabular}\vspace{-1mm}
\begin{center}Table 1: Experimental Parameters\end{center}
Theoretically and experimentally, we require that $\beta \in (0,1)$. Selection of $\beta\approx 1/4$ was done by consulting values considered in the original GTD experiments, and $\alpha = 8$ or $\alpha=10$ was based on the fact that larger step-sizes yield faster learning, so it is advantageous to use step-sizes $\alpha$ an order of magnitude larger than $\beta$.}

To obtain a benchmark policy for this task, we make use of trust region policy optimization \cite{schulman2015trust}. To evaluate value function estimates, we form an offline training set of state transitions and associated rewards by running this policy through consecutive episodes until we had one training trajectory of 5000 steps and then repeat this for 100 training trajectories to generate sample statistics. For ground truth, we generate one long trajectory of 10000 steps and randomly sample 2000 states from it. From each of these 2000 states, we apply the policy until episode termination and use the observed discounted return as  $\hat{V}_\pi(\bbx)$. Since our policy was deterministic, we only performed this procedure once per sampled state. For value function $V$, we define the percentage error metric:
%
$  \text{Percentage Error}(V) = ({1}/{2000}) \sum_{i=1}^{2000} | ({ V(\bbx_i) - \hat{V}_\pi(\bbx_i) })/{\hat{V}_\pi(\bbx_i) } |$
%

We compared PKGTD with a Gaussian kernel to two other techniques for policy evaluation that also use kernel-based value function representations: (1) Gaussian process temporal difference (GPTD) \cite{engel2003bayes}, and (2) gradient temporal difference (GTD) \cite{sutton2009convergent} using radial basis function (RBF) network features.  We fix a kernel bandwidth across all techniques, and select parameter values that yield the best results for each method (see Table 1).  For RBF feature generation, we use two fixed grids with different spacing.  The first was one for which GTD yielded a value function estimate with percentage error similar to that which we obtained using PKGTD (RBF-49), and the second was one which yielded a number of basis functions that was similar to what PKGTD selected (RBF-25). 

{Figure~\ref{fig:experiments} displays these results: on the left we show percentage error, a surrogate for Bellman evaluation error, versus training example, in which we observe that PKGTD yields fast and reliable learning. In the center figure, we show the number of points in the kernel dictionary (model size) over samples, which demonstrates that PKGTD only keeps past states needed to estimate the value function well, rather than statistically insignificant points. 
Overall, we note that GTD with fixed RBF features requires a much denser grid in order to reach the same Percentage Error as Algorithm \ref{alg:pkgtd}, and that adaptive instance selection results in both faster initial learning and smaller error.  Compared to GPTD, which chooses model points online according to a fixed linear-dependence criterion, PKGTD requires fewer model points and converges to a better estimate of the value function more quickly and stably.

Fig. 1 (right) displays a contour plot of the value function -- the x-axis denotes position, the y-axis denotes velocity, and bold black dots denote kernel dictionary elements, i.e., past visited states that are essential for representing the estimate of $V^\pi$. The contour plot suggests that low value states are when one has small velocity near position $-0.6$ (the bottom of the hill). Moreover, value progressively increases as speed increases away from the bottom of the hill towards the top. 

 }
 
\vspace{-3mm}
\section{Conclusion}\label{sec:discussion}
In this paper, we considered the problem of policy evaluation in infinite MDPs with value functions that belong to a RKHS.  To solve this problem, we extended recent SQG methods for compositional stochastic programming to a RKHS, and used the result, combined with greedy sparse subspace projection, in a new policy-evaluation procedure called PKGTD (Algorithm \ref{alg:pkgtd}).  Under diminishing step sizes, PKGTD solves Bellman's evaluation equation exactly under the hypothesis that its fixed point belongs to a RKHS (Theorem \ref{theorem:convergence_wp1}).  Under constant step sizes, we can further guarantee finite-memory approximations (Corollary \ref{corollary1}) that still exhibit mean convergence to a neighborhood of the optimal value function (Theorem \ref{theorem:constant_stepsize_convergence}).  In our Mountain Car experiments, PKGTD yields excellent sample efficiency and model complexity, and therefore holds promise for large state space problems common in robotics where fixed state-action space tiling may prove impractical.

\vspace{-2mm}

\bibliographystyle{IEEEtran}
\bibliography{bibliography}

\appendix

\subsection{Kernel Parameterization}\label{apx_parameterization}

\begin{algorithm}[t]
\caption{Destructive Kernel Orthogonal Matching Pursuit (KOMP) \hspace{-2mm}}
\begin{algorithmic}
\label{alg:komp}
\REQUIRE  function $\tilde{V}$ defined by dict. $\tbD \in \reals^{p \times \tilde{M}}$, coeffs. $\tbw \in \reals^{\tilde{M}}$, approx. budget  $\epsilon_t > 0$ \\
\STATE \textbf{initialize} $V=\tilde{V}$, dict. $\bbD = \tbD$ with indices $\ccalI$, model order $M=\tilde{M}$, coeffs.  $\bbw = \tbw$.
\WHILE{candidate dictionary is non-empty $\ccalI \neq \emptyset$}
{\FOR {$j=1,\dots,\tilde{M}$}
	\STATE Find minimal approx. error without dict. element $\bbd_j$ \vspace{-2mm}
	$$\gamma_j = \min_{\bbw_{\ccalI \setminus \{j\}}\in\reals^{{M}-1}} \|\tilde{V}(\cdot) - \sum_{k \in \ccalI \setminus \{j\}} w_k \kappa(\bbd_k, \cdot) \|_{\ccalH} \; .$$ \vspace{-5mm}
\ENDFOR}
	\STATE Find index minimizing error: $j^* = \argmin_{j \in \ccalI} \gamma_j$
	\INDSTATE{{\bf{if }} minimal error exceeds threshold $\gamma_{j^*}> \epsilon_t$}
	\INDSTATED{\bf stop} 
	\INDSTATE{\bf else} 
	
	\INDSTATED Prune dictionary $\bbD\leftarrow\bbD_{\ccalI \setminus \{j^*\}}$
	\INDSTATED Revise set $\ccalI \leftarrow \ccalI \setminus \{j^*\}$, model order ${M} \leftarrow {M}-1$. 	\vspace{-3mm}
	\INDSTATED Compute updated weights $\bbw$ defined by dict. $\bbD$
 	\vspace{-2mm}$$\bbw = \argmin_{\bbw \in \reals^{{M}}} \lVert \tilde{V}(\cdot) - \bbw^T\boldsymbol{\kappa}_{\bbD}(\cdot) \rVert_{\ccalH}$$\vspace{-6mm}
	\INDSTATE {\bf end}
\ENDWHILE	
\RETURN ${V},\bbD,\bbw$ of model order $M \leq \tilde{M}$ s.t. $\|V - \tilde{V} \|_{\ccalH}\leq \eps_t$
\end{algorithmic}
\end{algorithm}

{\bf Coefficient update} The update \eqref{eq:quasi_projected_fsgd}, for a fixed dictionary $\bbD_{t+1} \in \reals^{p\times M_{t+1}}$, may be expressed in terms of the parameter space of coefficients only.  To do so, first define the stochastic quasi-gradient update \emph{without projection}, given function ${V}_t$ parameterized by dictionary $\bbD_t$ and coefficients $\bbw_t$, as\vspace{-1mm}
\begin{align}\label{eq:sgd_tilde}
\tilde{V}_{t+1} =(1-\alpha_t\lambda){V}_t - \alpha_t (\gamma \kappa(\bby_t, \cdot) - \kappa(\bbx_t,\cdot))z_{t+1}  \; .
\end{align}
This update may be represented using dictionary and weights \vspace{-3mm}
\begin{align}\label{eq:param_tilde}
&\tbD_{t+1} = [\bbD_t \; , \; \bbx_t \; , \; \bby_t ] \nonumber\\
&\tbw_{t+1} = [ (1 - \alpha_t \lambda) \bbw_t \; , \; \alpha_t  z_{t+1} \; , \; - \alpha_t \gamma z_{t+1}]  \; ,
\end{align}
{Here we drop the time index for notational clarity but note that it can be easily inferred from the context. }
 $V_{t+1}$ denotes the projected {SQG} iterates [cf. \eqref{eq:quasi_projected_fsgd}] and whereas $\tilde{V}_{t+1}$ denotes the un-projected iterate [cf. \eqref{eq:sgd_tilde}] in Sec. \ref{subsec:polk_projection}. The later is parameterized by dictionary $\tbD_{t+1}$ and weights $\tbw_{t+1}$ \eqref{eq:param_tilde}. 

When the dictionary defining $V_{t+1}$ is assumed fixed, {we use} the Representer Theorem to rewrite \eqref{eq:projection_hat} as a kernel expansions, where the coefficients are the only free parameter:
\begin{align}\label{eq:proximal_hilbert_representer}
 &\argmin_{\bbw \in \reals^{{M}_{t+1}}} 
\frac{1}{2\eta_t} \Big\lVert \sum_{n=1}^{{M}_{t+1}} {w}_n\kappa(\bbd_n,\cdot)  - \sum_{m=1}^{\tilde{M}}\tilde{w}_m\kappa(\tbd_m,\cdot) \Big\rVert_{\ccalH}^2  \\
%
%
&\qquad\qquad:= \bbw_{t+1} \; . \nonumber
\end{align}
In \eqref{eq:proximal_hilbert_representer}, the first equality comes from expanding the square, and the second comes from defining the cross-kernel matrix $\bbK_{\bbD_{t+1},\tbD_{t+1}}$ whose $(n,m)^\text{th}$ entry is $\kappa(\bbd_n,\tbd_m)$.  Kernel matrices $\bbK_{\tbD_{t+1},\tbD_{t+1}}$ and $\bbK_{\bbD_{t+1},\bbD_{t+1}}$ are similarly defined. Here ${M}_{t+1}$ is the number of columns in $\bbD_{t+1}$, while $\tilde{M}_{t+1}=M_t + 2$ is that of in $\tbD_{t+1}$ [cf. \eqref{eq:param_tilde}]. 
Observe that $\tbD_{t+1}$ has $\tilde{M}_{t+1}=M_t + 2$ columns, which is the length of $\tbw_{t+1}$. For a fixed dictionary $\bbD_{t+1}$, the stochastic projection in \eqref{eq:projection_hat} is a least-squares problem on the coefficient vector, i.e.,

\begin{equation} \label{eq:hatparam_update}
\bbw_{t+1}=  \bbK_{\bbD_{t+1} \bbD_{t+1}}^{-1} \bbK_{\bbD_{t+1} \tbD_{t+1}} \tbw_{t+1} \;,
\end{equation}
 The explicit solution of \eqref{eq:proximal_hilbert_representer} may be obtained by noting that the last factor is independent of $\bbw$, and thus by computing gradients and solving for $\bbw_{t+1}$ we obtain \eqref{eq:hatparam_update}. Now we turn to {dictionary} selection $\bbD_{t+1}$ from trajectory $\{\bbx_u, \pi(\bbx_u), \bby_u\}_{u \leq t}$.

{\bf Dictionary Update} We select dictionary $\bbD_{t+1}$ {via} greedy compression, a topic studied in compressive sensing \cite{needell2008greedy}.  The function $\tilde{V}_{t+1} =(1-\alpha_t) {V}_t - \alpha_t (\gamma \kappa(\bby_t, \cdot) - \kappa(\bbx_t,\cdot))z_{t+1}$ defined by SQG method without projection \eqref{eq:sgd_tilde} is parameterized by dictionary $\tbD_{t+1}$ [cf. \eqref{eq:param_tilde}]. We form $\bbD_{t+1}$ by selecting a subset of $M_{t+1}$ columns from $\tbD_{t+1}$ that best approximate $\tilde{V}_{t+1}$ in terms of Hilbert norm error. {This specification may be met via \emph{kernel orthogonal matching pursuit} (KOMP) \cite{Vincent2002} with error tolerance $\epsilon_t$, which yields a dictionary $\bbD_{t+1}$ comprised of a subset of columns of $\tbD_{t+1}$}. We tune $\eps_t$ to ensure both descent (Lemma \ref{lemma:iterate_relations}\eqref{lemma:value_function_suboptimality}) and finite memory (Corollary \ref{corollary1}).

With respect to the KOMP procedure above, we specifically use a variant called destructive KOMP with pre-fitting (see \cite{Vincent2002}, Section 2.3). This flavor of KOMP takes as an input a candidate function $\tilde{V}$ of model order $\tilde{M}$ parameterized by its dictionary $\tbD\in\reals^{p\times\tilde{M}}$ and coefficients $\tbw\in\reals^{\tilde{M}}$. The method then approximates $\tilde{V}$ by $V\in \ccalH$ with a lower model order. Initially, the candidate is the original $V = \tilde{V} $ so that its dictionary is initialized with $\bbD=\tbD$, with coefficients  $\bbw=\tbw$. 
Then, we sequentially {and greedily} remove model points from initial dictionary $\tbD$ {until  threshold $\|V - \tilde{V} \|_{\ccalH} \leq \eps_t $ is violated.  The result is a sparse approximation $V$ of $\tilde{V}$.} {Moreover, we also assume that the $V_{t+1}$ output from KOMP has bounded Hilbert norm, which is often required in the analysis of stochastic optimization algorithms. This assumption can be explicitly enforced by adding a bounded norm constraint into the the optimization problem for finding the best set of bases in the matching pursuit algorithm, which attainable by thresholding the coefficient sequence during compression.}

This process is executed via destructive KOMP. At each stage, a single dictionary element $j$ of $\bbD$ is selected to be removed which contributes the least to the Hilbert-norm approximation error $\min_{V\in\ccalH_{\bbD\setminus \{j\}}}\|\tilde{V} - V \|_{\ccalH}$ of the original function $\tilde{V}$, when dictionary $\bbD$ is used. Since at each stage the kernel dictionary is fixed, this amounts to a computation involving weights $\bbw \in \reals^{M-1}$ only; that is, the error of removing dictionary point $\bbd_j$ is computed for each $j$ as 
 $\gamma_j =\min_{\bbw_{\ccalI \setminus \{j\}}\in\reals^{{M}-1}} \|\tilde{V}(\cdot) - \sum_{k \in \ccalI \setminus \{j\}} w_k \kappa(\bbd_k, \cdot) \|.$
 $\bbw_{\ccalI \setminus \{j\}}$ denotes the entries of $\bbw\in \reals^M$ restricted to the sub-vector associated with indices $\ccalI \setminus \{j\}$. Then, we define the dictionary element which contributes the least to the approximation error as $j^*=\argmin_j \gamma_j$. If the error associated with removing this kernel dictionary element exceeds the given approximation budget $\gamma_{j^*}>\eps_t$, the algorithm terminates. Otherwise, this dictionary element $\bbd_{j^*}$ is removed, the weights $\bbw$ are revised based on the pruned dictionary as $\bbw = \argmin_{\bbw \in \reals^{{M}}} \lVert \tilde{f}(\cdot) - \bbw^T\boldsymbol{\kappa}_{\bbD}(\cdot) \rVert_{\ccalH}$, and the process repeats as long as the current function approximation is defined by a nonempty dictionary. See Algorithm \ref{alg:komp} for a summary.


\subsection{Auxiliary Results and Technical Lemmas}\label{apx_assumptions}

Next we turn to establishing some technical results which are necessary precursors to the main stability results.

\begin{proposition}\label{prop1}
Given independent identical realizations $(\bbx_t, \pi(\bbx_t), \bby_t)$ of the random triple $(\bbx, \pi(\bbx), \bby)$, the difference between the projected stochastic functional quasi-gradient and the stochastic functional quasi-gradient of the instantaneous cost instantaneous risk defined by \eqref{eq:quasi_sg} and \eqref{eq:projected_quasi_sg}, respectively, is bounded for all $t$ as\vspace{-1mm}
\begin{align}\label{eq:prop1}
\!\!\| \tilde{\nabla}_V \!J(V_t, z_{t+1}; \!\bbx_t, \pi(\!\bbx_t\!),\! \bby_t\! )\! -\!\! \hat{\nabla}_V J\!(V_t, z_{t+1}; \!\bbx_t, \pi(\!\bbx_t\!), \!\bby_t\! )\!\|_{\ccalH}\!
 \leq \!\frac{\eps_t}{\alpha_t}
\end{align}
where $\alpha_t>0$ denotes the algorithm step-size and $\eps_t>0$ is the compression budget parameter of Algorithm \ref{alg:komp}.
\end{proposition}\vspace{-1mm}
\begin{myproof}
As in Proposition 6 of \cite{koppel2019parsimonious}, consider the square-Hilbert-norm difference of $ \tilde{\nabla}_V J(V_t, z_{t+1}; \bbx_t, \pi(\bbx_t), \bby_t ) $  [cf. \eqref{eq:quasi_sg}] and $\hat{\nabla}_V J(V_t, z_{t+1}; \bbx_t, \pi(\bbx_t), \bby_t ) $ [cf. \eqref{eq:projected_quasi_sg}] \vspace{-1mm}
\begin{align}\label{eq:norm_stoch_grad_diff}
& \| \tilde{\nabla}_V \!J(V_t, z_{t+1}; \!\bbx_t, \pi(\!\bbx_t\!),\! \bby_t\! )\! -\!\! \hat{\nabla}_V J\!(V_t, z_{t+1}; \!\bbx_t, \pi(\!\bbx_t\!), \!\bby_t\! )\!\|_{\ccalH}\!\nonumber \\
& \  = \Big\| \! \frac{1}{\alpha} \Big( V_{t}\! - \!\ccalP_{ \ccalH_{\bbD_{t+1}}} \!\!\Big[ \!V_t\! -\! \alpha_t \hat{\nabla}_V J(V_t, z_{t+1}; \bbx_t, \pi(\bbx_t), \bby_t )\!\Big]\!\Big)\! \nonumber\\
 &\qquad - \hat{\nabla}_V J(V_t, z_{t+1}; \bbx_t, \pi(\bbx_t), \bby_t ) \Big\|_{\ccalH}^2 
\end{align}
Multiply and divide $\hat{\nabla}_V J(V_t, z_{t+1}; \bbx_t, \pi(\bbx_t), \bby_t )$, the last term, by $\alpha_t$, and reorder terms to write
\begin{align}\label{eq:norm_stoch_grad_expand}
& \Bigg\| \frac{\left( V_{t}  - \alpha_t \hat{\nabla}_V J(V_t, z_{t+1}; \bbx_t, \pi(\bbx_t), \bby_t )\right)}{\alpha_t} \nonumber\\
 & \quad- \frac{\ccalP_{ \ccalH_{\bbD_{t+1}}} \Big[ V_t - \alpha_t \hat{\nabla}_V J(V_t, z_{t+1}; \bbx_t, \pi(\bbx_t), \bby_t )\Big]\Big)}{\alpha_t }
 \Bigg\|_{\ccalH}^2 \nonumber \\
  &=  \frac{1}{\alpha_t^2}\Big\|\! ( V_{t} \! - \!\alpha_t \hat{\nabla}_V J(V_t, z_{t+1}; \bbx_t, \pi(\bbx_t), \bby_t ) \nonumber\\
&\qquad - \ccalP_{ \ccalH_{\bbD_{t+1}}} \Big[ V_{t}  \!- \!\alpha_t \hat{\nabla}_V J(V_t, z_{t+1}; \bbx_t, \pi(\bbx_t), \bby_t )\Big]\Big) \Big\|_{\ccalH}^2 \nonumber \\
& = \frac{1}{\alpha_t^2}\|\tilde{V}_{t+1} - V_{t+1}\|_{\ccalH}^2 \leq  \frac{\epsilon_t^2}{\alpha_t^2}
\end{align}
  where we have pulled the nonnegative scalar $\alpha_t$ outside the norm on the second line and substituted the definition of $\tilde{V}_{t+1}$ and $V_{t+1}$ in \eqref{eq:quasi_fsgd} and \eqref{eq:quasi_projected_fsgd}, respectively, in the last one. These facts combined with the KOMP residual stopping criterion in Algorithm \ref{alg:komp} is $\lVert \tilde{V}_{t+1} - V_{t+1} \rVert_{\ccalH} \leq \epsilon_t$ applied to the last term on the right-hand side of \eqref{eq:norm_stoch_grad_expand} yields \eqref{eq:prop1}. 

\end{myproof}

\begin{lemma}\label{lemma:iterate_relations} 
Let Assumptions \ref{as:first} - \ref{as:mean_variance_frechet} hold true and consider the sequence of iterates defined by Algorithm \ref{alg:pkgtd}. Then:
 \begin{enumerate}
 
\item \label{lemma:value_function_difference}The conditional expectation of the Hilbert-norm difference of value functions at the next and current iteration satisfies the relationship
\begin{align}\label{eq:value_function_difference}
\!\!\!\mathbb{E}\!\left[ \|V_{t+1} \!-\! V_{t} \|_{\ccalH}^2 \given \ccalF_t \right] \!\leq \! 2\alpha_t^2 (G_{\delta}^2 G_{V}^2 \! +\!\lambda^2K^2) \!+\!2 \epsilon_t^2
\end{align}
\item \label{lemma:value_function_suboptimality}The conditional expectation of the Hilbert-norm difference of value functions at the next and current iteration satisfies the relationship
\begin{align}\label{eq:value_function_suboptimality}
\!\!\!\mathbb{E}\!\left[\|V_{t+1} \!-\!\! V^*\! \|_{\ccalH}^2\! \given\! \ccalF_t \right]\! \!&\leq 
 \!\!\left(\!\!1\! \!+\!  \frac{\alpha_t^2}{\beta_t} G_V^2\!\!\!\right)\!\! \| V_t\! -\! V^*\! \|_{\ccalH}^2\!  \\
 & \quad+\! 2 \eps_t \|V_t \!-\!\!V^{\!*} \!\|_{\ccalH}\! \!-\! 2 \alpha_t \! \!\left[  J\!(\!V_t\!)\! - \!\!J(\!V^{\!*}\!)  \right] \! \!\nonumber\\
&\quad +\! \alpha_t^2  \sigma_V^2 +  \beta_t  \mathbb{E}\left[(z_{t+1} - \bar{\delta}_t  )^2\! \given\! \ccalF_t \right]  .\nonumber
\end{align}
 \item \label{lemma:per_iterate_td_average} Define the expected value of the temporal difference given the state variable $\bbx$ and policy $\pi$ as $\bar{\delta}_t = \mathbb{E}[\delta_t \given \bbx_t, \pi(\bbx_t)]$. Then the evolution of the auxiliary sequence $z_t$ with respect to $\bar{\delta}_t $ satisfies
\begin{align} \label{eq:per_iterate_td_average}
\!\!\mathbb{E}\!\left[(z_{t+1} \!\!-\! \bar{\delta}_t)^2\! \given \!\ccalF_t \right]  
&\leq \!(\!1\!-\!\beta_t\!)(\!z_t\! - \!\bar{\delta}_{t-1}\!)^2 \!\!+\! \frac{L_V}{\beta_t} \!\| V_t \!- \!V_{t-1}\! \|_{\ccalH}^2 \nonumber\\
& \quad+ 2 \beta_t^2 \sigma_{\delta}^2 
\end{align}
\end{enumerate}
\end{lemma}
{\bf \noindent Proof of Lemma \ref{lemma:iterate_relations}\eqref{lemma:value_function_difference}}: Consider the Hilbert-norm difference of value functions at the next and current iteration, and use the definition of $V_{t+1}$ in \eqref{eq:quasi_projected_fsgd_projected}, i.e.,
\begin{align}\label{eq:value_function_difference_1}
\|V_{t+1} - V_{t} \|_{\ccalH}^2 &= \alpha_t^2 \|   \tilde{\nabla}_V J(V_t, z_{t+1}; \bbx_t, \pi(\bbx_t), \bby_t )   \|_{\ccalH}^2 \nonumber \\
& \leq 2\alpha_t^2 \|   \hat{\nabla}_V J(V_t, z_{t+1}; \bbx_t, \pi(\bbx_t), \bby_t )   \|_{\ccalH}^2  \nonumber \\
&\; + 2\alpha_t^2 \|   \hat{\nabla}_V J(V_t, z_{t+1}; \bbx_t, \pi(\bbx_t), \bby_t ) \nonumber\\
&\; - \tilde{\nabla}_V J(V_t, z_{t+1}; \bbx_t, \pi(\bbx_t), \bby_t )  \|_{\ccalH}^2 \; ,
\end{align}
where we add and subtract the functional stochastic quasi-gradient $\hat{\nabla}_V J(V_t, z_{t+1}; \bbx_t, \pi(\bbx_t), \bby_t )$ on the first line of \eqref{eq:value_function_difference_1} and {the fact that the square of a sum is less than the sum of squares, due to Cauchy-Schwartz, i.e., $(a+b)^2 \leq 2 a^2 + 2 b^2$ for any $a,b>0$}. Now, we may apply Proposition \ref{prop1} to the second term and compute the conditional expectation to obtain
\begin{align}\label{eq:value_function_difference_2}
&\mathbb{E}[\|V_{t+1} - V_{t} \|_{\ccalH}^2\given \ccalF_t ] \\
&\quad =2\alpha_t^2 \mathbb{E}[\|   \hat{\nabla}_V J(V_t, z_{t+1}; \bbx_t, \pi(\bbx_t), \bby_t )   \|_{\ccalH}^2 \given \ccalF_t]  + 2\epsilon_t^2 \; .
\end{align}
Use the Cauchy-Schwartz inequality together with Law of Total Expectation and the definition of the functional stochastic quasi-gradient \eqref{eq:quasi_sg} to upper-estimate \eqref{eq:value_function_difference_2} as
\begin{align}\label{eq:value_function_difference_3}
&\mathbb{E}[\|V_{t+1}\! - \!V_{t} \|_{\ccalH}^2\given \ccalF_t ]\!  \nonumber\\
&\quad\leq 2\alpha_t^2 \mathbb{E}\Big\{\! \|  \gamma \kappa(\bby_t, \!\cdot) \!-\! \kappa(\bbx_t,\!\cdot))\|^2_{\ccalH}  \nonumber \\
&\qquad \times \mathbb{E}[ z_{t+1}^2\! \given \bbx_t, \pi(\bbx_t)] \given \ccalF_t \!\Big\} \!\!+\!2\alpha_t^2\lambda \|V_t\|_{\ccalH}^2
 + 2\epsilon_t^2  \; ,
\end{align}
which together with equation \ref{eq:td_assumption} (Assumption \ref{as:mean_variance_td}) regarding fact that $z_{t+1}$ has a finite second conditional moment, yields
\begin{align}\label{eq:value_function_difference_4}
\mathbb{E}[\|V_{t+1}\! -\! V_{t} \|_{\ccalH}^2\!\given \!\ccalF_t] & \leq2\alpha_t^2 G_{\delta}^2 \mathbb{E}\Big[ \! \|  \gamma \kappa(\!\bby_t, \cdot)\! -\! \kappa(\bbx_t,\cdot\!)\!)\|^2_{\ccalH}\! \given\! \ccalF_t\! \Big] \nonumber\\
&\quad +2\alpha_t^2\lambda \|V_t\| + 2\epsilon_t^2 \nonumber \\
& \leq 2\alpha_t^2 (\!G_{\delta}^2 G_{V}^2 +\lambda^2 K^2) \!+ \!2\epsilon_t^2 \; ,
\end{align}
where we have also applied the fact that the functional gradient of the temporal difference $\gamma \kappa(\bby_t, \cdot) - \kappa(\bbx_t,\cdot))$ has a finite second conditional moment and the bound on the function sequence [cf. \eqref{eq:bounded_iterates}], allowing us to conclude \eqref{eq:value_function_difference}.  \hfill$\blacksquare$ \\ \\
\noindent {\bf Proof of Lemma \ref{lemma:iterate_relations}\eqref{lemma:value_function_suboptimality}}: This proof is a generalization of Lemma 3 in Appendix G.2 in the Supplementary Material of \cite{wang2017stochastic} to a function-valued stochastic quasi-gradient step combined with bias induced by the sparse subspace projections $\ccalP_{\ccalH_{\bbD_{t+1}}}[\cdot]$ in \eqref{eq:quasi_projected_fsgd}. Begin by considering the square-Hilbert norm sub-optimality of $V_{t+1}$, i.e.,
\begin{align}\label{eq:projected_fsgd_suboptimality1}
&\|V_{t+1} - V^* \|_{\ccalH}^2  \nonumber\\
&= \| V_t - \alpha_t \tilde{\nabla}_V J(V_t, z_{t+1}; \bbx_t, \pi(\bbx_t), \bby_t ) -V^*\|_{\ccalH}^2 \nonumber \\
&= \| V_t \!-\! V^* \!\|_{\ccalH}^2  \!-\! 2 \alpha_t \langle  \tilde{\nabla}_V J(V_t, z_{t+1}; \bbx_t, \pi(\!\bbx_t\!), \bby_t \!),\!\! V_t \!-\!\!V^*\! \rangle_{\ccalH} \nonumber \\
&\qquad + \alpha_t^2 \| \tilde{\nabla}_V J(V_t, z_{t+1}; \bbx_t, \pi(\bbx_t), \bby_t )\|_{\ccalH}^2 \; ,
\end{align}
where we use the reformulation of the projected functional stochastic quasi-gradient step defined in \eqref{eq:quasi_projected_fsgd_projected} for the first equality, and expand the square in the second. Now, adding and subtracting $\hat{\nabla}_V J(V_t, z_{t+1}; \bbx_t, \pi(\bbx_t), \bby_t )$ the (un-projected) functional stochastic quasi-gradient \eqref{eq:quasi_sg} yields
\begin{align}\label{eq:projected_fsgd_suboptimality2}
&\|V_{t+1} - V^* \|_{\ccalH}^2 \nonumber\\
&\quad=\! \| V_t \!- \!\!V^*\! \|_{\ccalH}^2 \!
\!-\! 2 \alpha_t \langle  \hat{\nabla}_V J(\!V_t, \!z_{t+1}; \bbx_t, \!\pi(\!\bbx_t\!), \bby_t\! ), \! V_t\!-\!V^*\! \rangle_{\ccalH} \nonumber \\
&\qquad\quad + 2 \alpha_t \langle  \hat{\nabla}_V J(V_t, \!z_{t+1}; \bbx_t, \pi(\bbx_t), \bby_t ) \nonumber\\
&\qquad\quad - \tilde{\nabla}_V J(V_t, z_{t+1}; \bbx_t, \pi(\bbx_t), \bby_t ), V_t \!-\!V^* \rangle_{\ccalH} \nonumber \\
&\qquad\quad+ \alpha_t^2 \| \tilde{\nabla}_V J(V_t, z_{t+1}; \bbx_t, \pi(\bbx_t), \bby_t )\|_{\ccalH}^2 \; .
\end{align}
Apply the Cauchy-Schwartz inequality to the third term on the right-hand side of \eqref{eq:projected_fsgd_suboptimality2} together with the bound on the difference between unprojected and projected stochastic quasi-gradients in Proposition \ref{prop1} to obtain
\begin{align}\label{eq:projected_fsgd_suboptimality3}
&\|V_{t+1} - V^* \|_{\ccalH}^2 \\
& = \| V_t \!-\! V^*\! \|_{\ccalH}^2
\! -\! 2 \alpha_t \langle  \hat{\nabla}_V J(\!V_t, z_{t+1}; \bbx_t, \pi(\!\bbx_t\!), \bby_t \!), \! V_t\!-\!V^* \rangle_{\ccalH} \nonumber \\
& \  + 2 \eps_t \|V_t -V^* \|_{\ccalH}+ \alpha_t^2 \| \tilde{\nabla}_V J(V_t, z_{t+1}; \bbx_t, \pi(\bbx_t), \bby_t )\|_{\ccalH}^2 \; .\nonumber
\end{align}
Now,  with $\bar{\delta}_t= \mathbb{E}[\delta_t \given \bbx_t, \pi(\bbx_t)]$, add and subtract $\hat{\nabla}_V J(V_t, \bar{\delta}_t; \bbx_t, \pi(\bbx_t), \bby_t )$, the stochastic quasi-gradient evaluated at $(V_t, \bar{\delta}_t)$ rather than $(V_t, z_{t+1})$, inside the inner-product term on the right-hand side of \eqref{eq:projected_fsgd_suboptimality3}, to write
\begin{align}\label{eq:projected_fsgd_suboptimality4}
&\|V_{t+1} - V^* \|_{\ccalH}^2 \nonumber\\
& = \| V_t - V^* \|_{\ccalH}^2
 - 2 \alpha_t \langle  \hat{\nabla}_V J(V_t, \delta_t; \bbx_t, \pi(\bbx_t), \bby_t ), V_t -V^* \rangle_{\ccalH} \nonumber\\
 &\quad + 2 \eps_t \|V_t -V^* \|_{\ccalH} 
+  2 \alpha_t \langle  (\gamma \kappa(\bby_t, \cdot) - \kappa(\bbx_t,\cdot))(\bar{\delta}_t - z_{t+1} ), \nonumber \\
&\quad \ V_t -V^* \rangle_{\ccalH} + \alpha_t^2 \| \tilde{\nabla}_V J(V_t, z_{t+1}; \bbx_t, \pi(\bbx_t), \bby_t )\|_{\ccalH}^2 \; ,
\end{align}
where we substitute in the definitions of $\hat{\nabla}_V J(V_t, \bar{\delta}_t; \bbx_t, \pi(\bbx_t), \bby_t )$ and $\hat{\nabla}_V J(V_t, z_{t+1}; \bbx_t, \pi(\bbx_t), \bby_t )$ [cf. \eqref{eq:stoch_grad}, \eqref{eq:quasi_sg}, respectively] in \eqref{eq:projected_fsgd_suboptimality4}, and cancel out the common regularization term $\lambda V_t$. We define the directional error associated with difference between the stochastic quasi-gradient and the stochastic gradient as 
\begin{equation}\label{eq:quasi_gradient_error}
v_t=2 \alpha_t \langle  (\gamma \kappa(\bby_t, \cdot) - \kappa(\bbx_t,\cdot))(\bar{\delta}_t - z_{t+1} ), V_t -V^* \rangle_{\ccalH}
\end{equation}
From here, compute the expectation conditional on $\ccalF_t$:
\begin{align}\label{eq:projected_fsgd_suboptimality5}
&\!\!\mathbb{E}\left[\|V_{t+1} - V^* \|_{\ccalH}^2 \given \ccalF_t \right] \nonumber\\
& \!\!\!\!= \!\|\! V_t \!-\! V^* \!\|_{\ccalH}^2\!\! -\! 2 \alpha_t \big\langle \mathbb{E}\!\left[ \!\hat{\nabla}_V \!J\!(\!V_t, \!\bar{\delta}_t; \!\bbx_t, \pi(\!\bbx_t\!), \bby_t\! )\!\given \!\ccalF_t \! \right]\!\!, \!\!V_t\! -\!V^* \!\rangle_{\ccalH}  \nonumber \\
&\quad  
 + 2 \eps_t \|V_t -V^* \|_{\ccalH}+  \mathbb{E}\left[ v_t\given \ccalF_t \right] \nonumber \\
 &\quad+ \alpha_t^2 \mathbb{E}\left[\| \tilde{\nabla}_V J(V_t, z_{t+1}; \bbx_t, \pi(\bbx_t), \bby_t )\|_{\ccalH}^2 \given \ccalF_t\right] \; .
\end{align}
Note that the compositional objective $J(V)$ is convex with respect to $V$, which allows us to write
\begin{align}\label{eq:convexity}
&\!\!\!\! \Big\langle \!\mathbb{E}\!\left[ \!\hat{\nabla}_V\! J\!(V_t, \bar{\delta}_t; \bbx_t, \pi(\!\bbx_t\!),\! \bby_t \!)\!\given \!\ccalF_t \!\right]\!, \!V_t \!-\!\!V^* \!\!\Big\rangle_{\!\!\ccalH}\! \!\geq \!J(V_t\!) \!-\! J(V^*\!)  .\!\!\!
\end{align}
Now, we may use Assumption \ref{as:mean_variance_td} [cf. \eqref{eq:projected_quasi_sg_variance}] regarding the finite conditional moments of the projected stochastic quasi-gradient to the last term in \eqref{eq:projected_fsgd_suboptimality5} so that it may be replaced by its upper-estimate, which together with \eqref{eq:convexity} simplifies to
\begin{align}\label{eq:projected_fsgd_suboptimality6}
\!\!\!\mathbb{E}\!\left[\|V_{t+1} \!\!-\! V^*\! \|_{\ccalH}^2\! \given\! \ccalF_t \right]  &\!=\! \| V_t\! -\! V^* \|_{\ccalH}^2
 \!-\! 2 \alpha_t \left[ J\!(V_t)\! -\! J(V^*) \right] \! \\
 &\quad\!+\! 2 \eps_t \|V_t \!-\!V^*\! \|_{\ccalH}\!+\! \alpha_t^2 \sigma_V^2 \! +\!  \mathbb{E}\left[ v_t\!\given\! \ccalF_t \right] .\nonumber
\end{align}
We need to analyze $v_t$, the directional error associated with using stochastic quasi-gradients rather than stochastic gradients. In doing so, we derive the fact that the sub-optimality $\|V_t - V^*\|$ is intrinsically coupled to the auxiliary sequence $(z_{t+1} - \bar{\delta}_t)$, the focus of Lemma \ref{lemma:iterate_relations}\eqref{lemma:per_iterate_td_average}. Proceed by applying Cauchy-Schwartz to \eqref{eq:quasi_gradient_error}, which allows us to write
\begin{align}\label{eq:quasi_gradient_error2}
\!\!\!\!v_t\leq2 \alpha_t \|\gamma \kappa(\bby_t, \cdot)\! -\! \kappa(\bbx_t,\cdot)\|_{\ccalH}^2 | z_{t+1}\!-\! \bar{\delta}_t  | \| V_t \!-\!V^* \|_{\ccalH}
\end{align}
Note that $2a b \leq \rho a^2 + b^2/\rho $ for $\rho, a,b>0$, which we apply to \eqref{eq:quasi_gradient_error2} with $a=   |z_{t+1}- \bar{\delta}_t | $, $b=\alpha_t \|\gamma \kappa(\bby_t, \cdot) - \kappa(\bbx_t,\cdot)\|_{\ccalH} \| V_t -V^* \|_{\ccalH}$, and $\rho = \beta_t$ so that \eqref{eq:quasi_gradient_error2} becomes
\begin{align}\label{eq:quasi_gradient_error3}
\!\!\!v_t\leq \beta_t  (z_{t+1}\!-\! \bar{\delta}_t\!)^2 \!+\!  \frac{\alpha_t^2}{\beta_t} \|\gamma \kappa(\bby_t, \!\cdot\!)\!\! -\! \kappa(\bbx_t,\!\cdot\!)\|_{\ccalH}^2 \|\! V_t \!-\!\!V^*\! \|_{\ccalH}^2  .
\end{align}
The conditional mean of $v_t$ [cf. \eqref{eq:quasi_gradient_error}], using \eqref{eq:quasi_gradient_error3}, is then
\begin{align}\label{eq:quasi_gradient_error_mean}
\mathbb{E}\left[v_t \given \ccalF_t \right] 
&\leq \beta_t  \mathbb{E}\left[(z_{t+1}- \bar{\delta}_t  )^2 \given \ccalF_t \right]  \\
& \quad \!+\!  \frac{\alpha_t^2}{\beta_t} \mathbb{E}\left[\|\gamma \kappa(\bby_t, \!\cdot\!) \!-\! \kappa(\bbx_t,\!\cdot\!)\|_{\ccalH}^2\!\given\! \ccalF_t\! \right]  \| V_t \!-\!V^*\! \|_{\ccalH}^2 \nonumber \\
& \leq \beta_t  \mathbb{E}\!\left[(z_{t+1}\!-\! \bar{\delta}_t\! )^2 \!\given \!\ccalF_t \!\right]  \!+ \! \frac{\alpha_t^2}{\beta_t} G_V^2 \| V_t -V^* \|_{\ccalH}^2 \; ,\nonumber
\end{align}
where we apply the finite variance property of the functional component of the stochastic gradient [cf. \eqref{eq:value_kernel_assumption}] for the final inequality \eqref{eq:quasi_gradient_error_mean}. Substitute \eqref{eq:quasi_gradient_error_mean} into \eqref{eq:projected_fsgd_suboptimality6} and gather terms:
\begin{align}\label{eq:projected_fsgd_suboptimality_final}
&\mathbb{E}\left[\|V_{t+1} - V^* \|_{\ccalH}^2 \given \ccalF_t \right] \\
&\quad \leq \left(1 +  \frac{\alpha_t^2}{\beta_t} G_V^2\right) \| V_t - V^* \|_{\ccalH}^2 + 2 \eps_t \|V_t -V^* \|_{\ccalH} \nonumber\\
&\qquad-\! 2 \alpha_t \!\left[ J\!(V_t) \!- \!J(V^*\!) \right]\! +\! \alpha_t^2 \sigma_V^2\!+\!  \beta_t  \mathbb{E}\left[(z_{t+1}\!\!-\! \bar{\delta}_t \! )^2\! \given\! \ccalF_t \right]  ,\nonumber
\end{align}
which is as stated in Lemma \ref{lemma:iterate_relations}\eqref{lemma:value_function_suboptimality}.
 \hfill$\blacksquare$ \\ 

\noindent {\bf Proof of Lemma \ref{lemma:iterate_relations}\eqref{lemma:per_iterate_td_average}}: This proof is an adaptation of Lemma 2 in Appendix G.1 in the Supplementary Material of \cite{wang2017stochastic} to the recursively averaged temporal difference sequence $z_t$ defined in \eqref{eq:td_average}. Begin by defining the scalar quantity $e_t$ as the difference of mean temporal differences scaled by the forgetting factor $\beta_t$, i.e.
$e_t = (1-\beta_t)(\bar{\delta}_t - \bar{\delta}_{t-1})$. Then we consider the difference of the evolution of the auxiliary variable $z_{t+1}$ with respect to the conditional mean temporal difference $\bar{\delta}_t$, plus the difference of mean temporal differences: 
\begin{align}\label{eq:auxiliary_sequence_difference}
&z_{t+1} - \bar{\delta}_t + e_t \nonumber \\
&= (\!1\!-\!\beta_t) z_t\! +\! \beta_t \delta_t \!-\! [(1\!-\!\beta_t)\bar{\delta}_t \!+\! \beta_t \bar{\delta}_t]  \!+\! (1\!-\!\beta_t)(\bar{\delta}_t \!-\! \bar{\delta}_{t-1}) \nonumber\\
&=(1-\beta_t)\left( z_t - \bar{\delta}_{t-1} \right)  + \beta_t (\delta_t - \bar{\delta}_t)
 \end{align}
where we make use of the definition of $z_{t+1}$ in \eqref{eq:td_average}, the fact that $\bar{\delta}_t =[(1-\beta_t)\bar{\delta}_t + \beta_t \bar{\delta}_t]$, and the definition of $e_t$ on the first line of \eqref{eq:auxiliary_sequence_difference}, and in the second we gather terms with respect to coefficients $(1-\beta_t)$ and $\beta_t$, and cancel the redundant $ \bar{\delta}_t$ term. Now, consider the square of the expression \eqref{eq:auxiliary_sequence_difference}, using it's simplification on the right-hand side of the preceding expression
\begin{align}\label{eq:auxiliary_sequence_difference_square}
(z_{t+1}\! -\! \bar{\delta}_t + e_t )^2 & = [ (1\!-\!\beta_t)\left( z_t \!- \!\bar{\delta}_{t-1} \right)  \!+\! \beta_t (\delta_t - \bar{\delta}_t)]^2  \\
& \quad= (1-\beta_t)^2\left( z_t - \bar{\delta}_{t-1} \right)^2  + \beta_t^2 (\delta_t - \bar{\delta}_t)^2  \nonumber \\
&\qquad+ 2(1-\beta_t)\beta_t \left( z_t - \bar{\delta}_{t-1} \right)(\delta_t - \bar{\delta}_t)  .\nonumber
 \end{align}
where we expand the square to obtain the second line in the previous expression. Now, compute the expectation of \eqref{eq:auxiliary_sequence_difference_square} conditional on the filtration $\ccalF_t$, which yields
\begin{align}\label{eq:auxiliary_sequence_mean_square}
&\mathbb{E}[(z_{t+1} - \bar{\delta}_t + e_t )^2\given \ccalF_t ] \nonumber\\
&\quad = (1-\beta_t)^2\left( z_t - \bar{\delta}_{t-1} \right)^2  + \beta_t^2 \mathbb{E}[(\delta_t - \bar{\delta}_t)^2\given \ccalF_t ] \nonumber \\
&\qquad + 2(1-\beta_t)\beta_t \left( z_t - \bar{\delta}_{t-1} \right)\mathbb{E}[(\delta_t - \bar{\delta}_t) \given \ccalF_t ]  \; .
 \end{align}
Now we apply the assumption [cf. \eqref{eq:td_assumption}]  that the fact that the temporal difference $\delta_t$ is an unbiased estimator for its conditional mean $\bar{\delta}_t$ (so that the last term in the previous expression is null), with finite variance $\mathbb{E}[(\delta_t - \bar{\delta}_t)^2\given \ccalF_t ] \leq \sigma_{\delta}^2$ (Assumption \ref{as:mean_variance_td}), to write
\begin{align}\label{eq:auxiliary_sequence_mean_square2}
\!\!\!\mathbb{E}[(\!z_{t+1} \!-\! \bar{\delta}_t + e_t )^2\!\given\! \ccalF_t ] \!= \!(1\!-\!\beta_t)^2\!\left( \!z_t \!-\! \bar{\delta}_{t-1} \right)^2  \!\!+\! \beta_t^2 \sigma_{\delta}^2  .
 \end{align}
We may use the relationship in \eqref{eq:auxiliary_sequence_mean_square2} to obtain an upper estimate on the conditional mean square of $z_{t+1} - \bar{\delta}_t $ by using the inequality $\| a + b \|^2 \leq (1+\rho )\|a\|^2 + (1 + 1/\rho)\|b\|^2$ which holds for any $\rho > 0$: set $a= z_{t+1} - \bar{\delta}_t + e_t$, $b=-e_t$, and $\rho=\beta_t$. Therefore, we obtain
\begin{align}\label{eq:aux_sequence_generalized_triangle_ineq}
(z_{t+1} \!-\! \bar{\delta}_t  \!)^2 \!\leq\! (1\!+\! \beta_t)(z_{t+1} \!-\! \bar{\delta}_t \!+\! e_t )^2 \!+\! \left(\!\!1\! +\! \frac{1}{\beta_t} \!\right)\!e_t^2 \; .
 \end{align}
Now, we use the expected value of \eqref{eq:aux_sequence_generalized_triangle_ineq} in lieu of \eqref{eq:auxiliary_sequence_mean_square2}, while gaining a multiplicative factor of $(1+\beta_t)$ on the right-hand side of \eqref{eq:auxiliary_sequence_mean_square2} plus the error term $(1 + 1/\beta_t )e_t$, yielding
\begin{align}
\label{eq:auxiliary_sequence_mean_square3}
&\mathbb{E}[(z_{t+1} - \bar{\delta}_t  )^2\given \ccalF_t ] \\
& =(\!1\!+\!\beta_t)\! \Big[\!(\!1\!-\!\beta_t)^2\!\!\left( z_t - \bar{\delta}_{t-1} \right)^2  + \beta_t^2 \sigma_{\delta}^2\Big] + \left(\frac{1+ \beta_t}{\beta_t} \right)e_t^2  .\nonumber
 \end{align}
Apply the fact that $(1-\beta_t^2)(1-\beta_t) \leq (1-\beta_t)$ to the first term in \eqref{eq:auxiliary_sequence_mean_square3} and  $(1+\beta_t)\beta_t^2 \leq 2 \beta_t^2$ to the second (since $\beta_t \in (0,1)$) to simplify \eqref{eq:auxiliary_sequence_mean_square3} as
\begin{align}\label{eq:auxiliary_sequence_mean_square4}
&\mathbb{E}[(z_{t+1} - \bar{\delta}_t  )^2\given \ccalF_t ]\\
& \qquad=(1-\beta_t)\left( z_t - \bar{\delta}_{t-1} \right)^2  + 2 \beta_t^2 \sigma_{\delta}^2 
 + \left(\frac{1+ \beta_t}{\beta_t} \right)e_t^2   .\nonumber
 \end{align}
Now we analyze the term involving $e_t$, which represents the difference of mean temporal differences. By definition,
\begin{align}\label{eq:mean_temporal_difference_difference}
\!\!\!| e_t | = (1\!-\!\beta_t)|(\bar{\delta}_t \!-\! \bar{\delta}_{t-1})| 
\! \leq \!(1\!-\!\beta_t) L_{V} \| V_t \!-\! V_{t-1} \|_{\ccalH}
\end{align}
where we apply the Lipschitz continuity of the conditional average temporal difference $ \bar{\delta}_t = \mathbb{E}_{\bby_t} [ r(\bbx_t, \pi(\bbx_t), \bby_t) + \gamma V(\bby_t) - V(\bbx_t) \given \bbx_t, \pi(\bbx_t)] $ with respect to the value function [cf. \eqref{eq:td_lipschitz}] stated in Assumption \ref{as:mean_variance_frechet}. Substitute the right-hand side of \eqref{eq:mean_temporal_difference_difference} into \eqref{eq:auxiliary_sequence_mean_square4}, and simplify the expression in the last term as $(1-\beta_t^2)/\beta_t \leq 1/\beta_t$ to conclude \eqref{eq:per_iterate_td_average}.  \hfill$\blacksquare$ \\ \\
\begin{lemma}\label{lemma:coupled_supermartingale}(Coupled Supermartingale Theorem \cite{wang2013incremental}[Lemma 6]) Let $\{\xi_k\}$, $\{\zeta_k\}$, $\{u_k\}$, $\{\bar{u}_k\}$, $\{\eta_k\}$, $\{\theta_k\}$, $\{\varepsilon_k\}$, $\{\mu_k\}$, $\{\nu_k\}$ be sequences of nonnegative random variables such that
\begin{align}
\mathbb{E}[\xi_{k+1} \given \ccalG_k ] &\leq (1+\eta_k) \xi_k - u_k + c \theta_k \zeta_k + \mu_k \; , \label{eq:supermartingale_value_function}  \\
\mathbb{E}[\zeta_{k+1} \given \ccalG_k ] &\leq (1-\theta_k) \zeta_k - \bar{u}_k + \varepsilon_k \xi_k + \nu_k \; ,
\label{eq:supermartingale_temporal_difference}
\end{align}
where $\ccalG_k = \{\xi_s, \zeta_s, u_s, \bar{u}_s, \eta_s, \theta_s, \varepsilon_s, \mu_s, \nu_s \}_{s=0}^k$ is the filtration, and $c>0$ is a scalar. Assume the following conditions:
\begin{align} \label{eq:summability}
&\sum_{k=0}^\infty \eta_k < \infty\; , \ \quad \sum_{k=0}^\infty \varepsilon_k < \infty \; , \nonumber\\
&\ \sum_{k=0}^\infty \mu_k < \infty \; ,\ \sum_{k=0}^\infty \nu_k < \infty \;,
\end{align}
almost surely. Then $\xi_k$ and $\zeta_k$ converge almost surely to two nonnegative random variables, and we may conclude that
\begin{equation} \label{eq:summability_convergence}
\sum_{k=0}^\infty u_k < \infty\; , \ \quad \sum_{k=0}^\infty \bar{u}_k < \infty \; , \ \sum_{k=0}^\infty \theta_k \zeta_k < \infty \; \text{ a. s. } %
\end{equation}
%
\end{lemma}
We use Lemma \ref{lemma:coupled_supermartingale} to establish convergence w.p.$1$ of Algorithm \ref{alg:pkgtd} through the expressions derived in Lemma \ref{lemma:iterate_relations}.

\subsection{Proof of Theorem \ref{theorem:convergence_wp1}}\label{apx_theorem_convergence_wp1}

We use the relations established in Lemma \ref{lemma:iterate_relations} to construct a coupled supermartingale of the form in Lemma \ref{lemma:coupled_supermartingale} as follows. First, consider the expression \eqref{eq:value_function_suboptimality} for the value function sub-optimality, using approximation budget $\eps_t = \alpha_t^2$ and the fact that the value function is bounded in Hilbert norm [cf. \eqref{eq:bounded_iterates}] to obtain $ \|V_t -V^* \|_{\ccalH}\leq 2K$ :
\begin{align}\label{eq:value_function_suboptimality_supermartingale}
&\mathbb{E}\left[\|V_{t+1} - V^* \|_{\ccalH}^2 \given \ccalF_t \right] \nonumber\\
&\quad \leq \left(1\! +\!  \frac{\alpha_t^2}{\beta_t} G_V^2\right) \| V_t \!-\! V^* \|_{\ccalH}^2 \! -\! 2 \alpha_t \left[ J(V_t)\! - \!J(V^*) \right]  \nonumber \\
&\qquad+ \alpha_t^2 ( \sigma_V^2 + 4K)   +  \beta_t  \mathbb{E}\left[(z_{t+1} - \bar{\delta}_t  )^2 \given \ccalF_t \right] \; .
\end{align}
and then substitute \eqref{eq:per_iterate_td_average} regarding the evolution of $z_t$ with respect to its conditional expectation into \eqref{eq:value_function_suboptimality_supermartingale} to obtain :

\begin{align}\label{eq:value_function_suboptimality_supermartingale2}
&\mathbb{E}\left[\|V_{t+1} - V^* \|_{\ccalH}^2 \given \ccalF_t \right] \nonumber\\
&\quad \leq \left(1 +  \frac{\alpha_t^2}{\beta_t} G_V^2\right) \| V_t - V^* \|_{\ccalH}^2 - 2 \alpha_t \left[ J(V_t) - J(V^*) \right]   \nonumber\\
&\qquad+ \alpha_t^2 ( \sigma_V^2 + 4K)
 \! +\!  \beta_t  (1-\beta_t)(z_t - \bar{\delta}_{t-1})^2\! \nonumber\\
 &\qquad +\! L_V \| V_t - V_{t-1} \|_{\ccalH}^2 \!+\! 2 \beta_t^3 \sigma_{\delta}^2 \; .
\end{align}
Assume $\beta_t\in (0,1)$. Thus, the right-hand side of \eqref{eq:value_function_suboptimality_supermartingale2} yields

\begin{align}\label{eq:value_function_suboptimality_supermartingale3}
&\mathbb{E}\left[\|V_{t+1} - V^* \|_{\ccalH}^2 \given \ccalF_t \right] \nonumber\\
&\quad \leq \left(1 +  \frac{\alpha_t^2}{\beta_t} G_V^2\right) \| V_t - V^* \|_{\ccalH}^2 - 2 \alpha_t \left[ J(V_t) - J(V^*) \right]   \nonumber\\
 &\qquad
  +  \beta_t (z_t - \bar{\delta}_{t-1})^2 + \alpha_t^2 ( \sigma_V^2 + 4K) \nonumber\\
  &\qquad +  L_V \| V_t - V_{t-1} \|_{\ccalH}^2 + 2 \beta_t^2 \sigma_{\delta}^2 \; .
\end{align}
We may identify \eqref{eq:value_function_suboptimality_supermartingale3} with the first supermartingale relationship in Lemma \ref{lemma:coupled_supermartingale} [cf. \eqref{eq:supermartingale_value_function}] via the identifications
\begin{align}\label{eq:supermartingale_identification1}
&\xi_t = \|V_{t} - V^* \|_{\ccalH}^2 \; , \eta_t \!=\!  \frac{\alpha_t^2}{\beta_t} G_V^2\; 
,u_t \!=\! 2 \alpha_t \!\left[ J(V_t) \!-\! J(V^*\!) \right]  ,  \nonumber \\
&c=1\; ,\qquad\qquad\quad \zeta_t = (z_t - \bar{\delta}_{t-1})^2 \; ,\quad \theta_t = \beta_t \; ,\nonumber \\
& \mu_t =  \alpha_t^2 ( \sigma_V^2 + 4K) +  L_V \| V_t - V_{t-1} \|_{\ccalH}^2 + 2 \beta_t^2 \sigma_{\delta}^2\; ,
\end{align}
where $u_t \geq 0$ by the definition of the optimal objective $J(V^*)$.  To establish the summability of $\mu_t$, 
consider Lemma \ref{lemma:iterate_relations}\eqref{lemma:value_function_difference}, which establishes that t $\|V_t - V_{t-1}\|_{\ccalH}\leq\ccalO(\alpha_{t-1}^2)$. {Since $\sum_t \alpha_t^2< \infty $ [cf. \eqref{eq:step_size_condition}], we can sum both sides over all $t$  to conclude the series is finite in conditional expectation:
\begin{equation}\label{eq:summability_value_func_diff}
\sum \mathbb{E}[ \|V_t - V_{t-1}\|_{\ccalH} \given\ccalF_t ] \leq \alpha_{t-1}^2  < \infty .
\end{equation}
Now, rewrite \eqref{eq:summability_value_func_diff} with total expectation by selecting $\ccalF_0$. Note that since the individual terms $\| V_t - V_{t-1} \|_{\ccalH}^2$ are finite due to the stipulation that the output of KOMP yields finite Hilbert norm value functions, and non-negative by the definition of a norm, we can interchange the expectation (integral) and sum using the Monotone Convergence Theorem to conclude that 
\begin{equation}\label{eq:summability_value_func_diff_MCT}
 \mathbb{E}\left[ \sum \|V_t - V_{t-1}\|_{\ccalH}  \right] <  \infty .
\end{equation}
Thus,  $\sum_{t=0}^\infty \| V_t - V_{t-1} \|_{\ccalH}^2 < \infty$ w.p.$1$, implying $\sum_t \mu_t < \infty$.}

 Now, let's connect the evolution of the auxiliary temporal difference sequence $z_t$ \eqref{eq:td_average} in Lemma \ref{lemma:iterate_relations}\eqref{lemma:per_iterate_td_average}. In particular, \eqref{eq:per_iterate_td_average} is related to \eqref{eq:supermartingale_temporal_difference} via the identifications:
\begin{align}\label{eq:supermartingale_identification2}
\bar{u}_t = 0\; , \varepsilon_t = 0   \; , 
\nu_t= \frac{L_V}{\beta_t} \| V_t - V_{t-1} \|_{\ccalH}^2 + 2 \beta_t^2 \sigma_{\delta}^2 \; , 
\end{align}
with $\zeta_t =(z_t - \bar{\delta}_{t-1})^2 $ and $\theta_t = \beta_t$ as in \eqref{eq:supermartingale_identification1}. The summability of $\nu_t$ follows the following logic: consider the expression $ \| V_t - V_{t-1} \|_{\ccalH}^2/\beta_t$ of order $\ccalO(\alpha_t^2 /\beta_t)$ in conditional expectation by Lemma \ref{lemma:iterate_relations}\eqref{lemma:value_function_difference}. Sum the resulting conditional expectation for all $t$, which by the summability of the sequence $\sum_t \alpha_t^2/\beta_t<\infty$ is finite. Therefore, $\sum_t   \| V_t - V_{t-1} \|_{\ccalH}^2 /\beta_t<\infty $ almost surely. 

Together with the conditions on the step-size sequences $\alpha_t$ and $\beta_t$ \eqref{eq:step_size_condition}, the summability conditions \eqref{eq:summability} of Lemma \ref{lemma:coupled_supermartingale}, the Coupled Supermartingale Theorem, are satisfied, which implies that $\xi_t = \|V_{t} - V^* \|_{\ccalH}^2$ and $\zeta_t= (z_t - \bar{\delta}_{t-1})^2$ converge to two nonnegative random variables almost surely, and that:
\begin{align}\label{eq:supermartingale_conclusion}
\sum_t  \!\alpha_t \!\left[ J(V_t)\! -\! J(V^*) \right] \!<\! \infty \; ,\quad \sum_t\! \beta_t (z_{t+1}\! -\! \bar{\delta}_{t})^2 \!<\! \infty \; ,
\end{align}
almost surely. The non-summability of the step-size sequences $\alpha_t$ and $\beta_t$ \eqref{eq:step_size_condition} allows us to conclude that:
\begin{align}\label{eq:supermartingale_conclusion2}
\liminf_{t\rightarrow \infty}   J(V_t)  = J(V^*)   \; , \ \
\liminf_{t\rightarrow \infty} (z_{t+1} - \bar{\delta}_{t})^2 = 0\; ,
\end{align}
almost surely, and that $\|V_t - V^*\|^2_{\ccalH}$ converges to a nonnegative random variable with probability $1$, as does $(z_{t+1} - \bar{\delta}_{t})^2$.
We proceed to show that the entire sequence must converge. The rest of this proof is analogous to \cite{wang2017stochastic}, but is repeated here
for completeness.
Let $\Omega_{V^*}$ be the collection of sample paths such that  $\Omega_{V^*}=\{ \bby: \lim_t \| V_t (\bby) - V^*\| \text{ exists } \}$. Here we use the notation not that the value function is evaluated at state $\bby$ but instead is a function of random variable $\bby$. We just established above that $\mathbb{P}(\Omega_{V^*}) = 1$ for any $V^*\in\ccalH$. To prove that any limiting value function is optimal, we need to establish that $\cap_{V^*\in \ccalH} \Omega_{V^*}$ is measurable and $\mathbb{P}(\cap_{V^*\in \ccalH} \Omega_{V^*})=1$. 

To do so, note that since $J$ is convex, the set of minimizers of $J$, denoted as $\ccalH^* \subset \ccalH$, is separable, and has a countably dense subset $\ccalH_{Q}^*$. Thus the probability of divergence for some $V^*\in\ccalH_Q^*$ is the probability of a union of countably many sets, each having null probability. Therefore, we may write
\begin{equation}\label{eq:prob_divergence}
\mathbb{P}\left(\cap_{\ccalH_Q^*} \Omega_{V^*} \right) = 1- \mathbb{P}\left(\cup_{\ccalH^*_Q} \Omega_{V^*}^c \right) \geq 1 -  \sum_{V^*\in\ccalH^*_Q} \mathbb{P}\left( \Omega_{V^*}^c \right) = 1
\end{equation}
by simple application of De Morgan's Law and Boole's inequality. Then consider any $\tilde{V} \in \ccalH^\star$ which is the limit of a sequence of optimal value functions $\{\tilde{V}_k\}_{k=1}^\infty \subset\ccalH^\star$. We can prove that $\|\tilde{V}_t(\bby) - \tilde{V}\|$ is convergent provided that $\|\tilde{V}_t(\bby)- \tilde{V}_k\| $ is convergent for all $k$. Note that
\begin{align}\label{eq:triangle_limit_points}
\|V_t(\bby) - \tilde{V}_k \|_{\ccalH} &- \|\tilde{V}_k - \tilde{V}\|_{\ccalH} \nonumber \\
&\leq \|V_t(\bby) - \tilde{V} \|_{\ccalH} \nonumber  \\
& \leq  \|V_t(\omega) - \tilde{V}_k \|_{\ccalH} + \|\tilde{V}_k - \tilde{V} \|_{\ccalH}\; . 
\end{align}
Since $\|{V}_t(\bby) - \tilde{V}_k\|_{\ccalH}$ has a limit, take $t\rightarrow \infty$ in \eqref{eq:triangle_limit_points}, yielding:
\begin{align}\label{eq:triangle_limit_points_limit}
\lim_{t\rightarrow \infty} \|V_t(\bby) - \tilde{V}_k \|_{\ccalH}& - \|\tilde{V}_k - \tilde{V}\|_{\ccalH} 
\leq \liminf_{t\rightarrow \infty} \|V_t(\bby) - \tilde{V} \|_{\ccalH}  \nonumber \\
&\leq \limsup_{t\rightarrow \infty} \|V_t(\bby) - \tilde{V} \|_{\ccalH}  \\
&\leq \lim_{t\rightarrow \infty}  \|V_t(\bby) \!-\!\! \tilde{V}_k \|_{\ccalH} + \|\tilde{V}_k - \tilde{V} \|_{\ccalH}\; ,\nonumber
\end{align}
which, by subtracting $\liminf_{t\rightarrow \infty} \|V_t(\bby) - \tilde{V} \|_{\ccalH}   $ from both sides in \eqref{eq:triangle_limit_points_limit}, cancelling the common $\lim_{t\rightarrow \infty}  \|V_t(\bby) - \tilde{V}_k \|_{\ccalH} $, and combining terms, allows us to write
\begin{align}\label{eq:limsup_liminf_diff}
\limsup_{t\rightarrow \infty} \! \|V_t(\bby) \!- \!\tilde{V} \|_{\ccalH} \!-\!\liminf_{t\rightarrow \infty} \|V_t(\bby) \!- \!\tilde{V} \|_{\ccalH} \! \leq 2 \|\tilde{V}_k \!- \! \tilde{V} \|_{\ccalH}  .
\end{align}
Take $k \rightarrow \infty$ in \eqref{eq:limsup_liminf_diff}, for which $ \|\tilde{V}_t - \tilde{V} \|_{\ccalH} \rightarrow 0$, hence
\begin{align}\label{eq:limsup_liminf_equal}
\limsup_{t\rightarrow \infty} \|V_t(\bby) - \tilde{V} \|_{\ccalH} =\liminf_{t\rightarrow \infty} \|V_t(\bby) - \tilde{V} \|_{\ccalH}\; ,
\end{align}
and therefore $\|V_t(\bby) - \tilde{V} \|_{\ccalH} $ has a limit, so $\bby \in \Omega_{\tilde{V}^*}$, and therefore $\cap_{\ccalH^*_{Q}} \Omega_{V^*} \subset \Omega_{\tilde{V}}$. Consequently, $ \mathbb{P}\left(\cap_{\ccalH^*_{Q}} \Omega_{V^*} \right) = 1$. As a result, we have $ \left(\cap_{\ccalH^*} \Omega_{V^*} \right)^c \subset \left(\cap_{\ccalH^*_{Q}} \Omega_{V^*} \right)^c $, both of which are measurable and have null probability:  $\mathbb{P}\left((\cap_{\ccalH^*} \Omega_{V^*} )^c\right) \leq \mathbb{P}\left((\cap_{\ccalH^*_{Q}} \Omega_{V^*} )^c\right)=0$. Thus, $(\cap_{\ccalH^*} \Omega_{V^*} )$ is measurable and occurs with probability $1$. Put another way, $\|V_t - \tilde{V} \|_{\ccalH}$ is convergent for all $\tilde{V}\in\ccalH^*$ with probability $1$.

Now, we can use this fact together with \eqref{eq:supermartingale_conclusion2}, namely, $\liminf_{t\rightarrow \infty}   J(V_t)  = J(V^*)$, to establish that $V_t$ converges to the minimizer of $J(V)$ a.s. To do so, let $V^*\in\ccalH^*$ the set of optimizers of $J$. Since $\|V_t(\bby) - V^*\|_{\ccalH}$ converges, it is bounded. Then, $\{V_t(\bby) \}$ must have a limit point $\tilde{V}$ being an optimal solution, $J(\tilde{V}) = J^*$ with $\tilde{V}\in\ccalH^*$, by the continuity of $J$. Since $\omega \in\cap{\ccalH^*} \Omega_{V^*} \subset \Omega_{\tilde{V}}$, $\{\|V_t(\bby) - \tilde{V}\|_{\ccalH}\}$ is a convergent sequence whose limit is null. Thus, $\|V_t(\bby) - \tilde{V}\|_{\ccalH}\rightarrow 0$, so $V_t(\bby) \rightarrow \tilde{V}$ on this sample path. $\tilde{V}$ is a random variable dependent on the sample path, parameterized by $\bby$. The set of all such sample paths has prob. $1$, so that $V_t$ converges to a random point in $\ccalH^*$.
\hfill$\blacksquare$ \\


\subsection{Proof of Theorem \ref{theorem:constant_stepsize_convergence}}\label{apx_theorem_constant_stepsize_convergence}
Before analyzing the mean convergence behavior of the value function, we consider the mean sub-optimality of the auxiliary variable $z_t$ with respect to the conditional mean of the temporal difference $\bar{\delta}_t$. To do so, compute the total expectation of Lemma \ref{lemma:iterate_relations}\eqref{lemma:per_iterate_td_average}, stated as
\begin{align} \label{eq:per_iterate_td_average2}
&\mathbb{E}\left[(z_{t+1} - \bar{\delta}_t)^2  \right] \\
& \  \!\leq \!(1\!-\!\beta)\E{(z_t \!-\! \bar{\delta}_{t-1})^2} \!+\! \frac{L_V}{\beta} \E{ \| V_t\! -\! V_{t-1} \|_{\ccalH}^2}\!\!+\! 2 \beta^2 \sigma_{\delta}^2 \; ,\nonumber
\end{align}
where we have substituted in constant learning rate $\beta_t = \beta$ in \eqref{eq:per_iterate_td_average2}. The total expectation of Lemma \ref{lemma:iterate_relations}\eqref{lemma:value_function_difference} regarding $ \| V_t - V_{t-1} \|_{\ccalH}^2$, the difference of value functions in Hilbert-norm, may be substituted into \eqref{eq:per_iterate_td_average2}, with constant step-size $\alpha_t=\alpha$ and compression budgets $\eps_t = \eps$ to obtain
\begin{align} \label{eq:per_iterate_td_average3}
&\mathbb{E}\left[(z_{t+1} - \bar{\delta}_t)^2  \right] \nonumber\\
&\quad \leq (1-\beta)\E{(z_t - \bar{\delta}_{t-1})^2} \nonumber\\
&\qquad + \frac{2 L_V}{\beta}\left[ \alpha^2 (G_{\delta}^2 G_{V}^2 +\lambda^2K^2) + \epsilon^2\right]
 + 2 \beta^2 \sigma_{\delta}^2 \; ,
\end{align}
Observe that \eqref{eq:per_iterate_td_average3} gives a relationship between the sequence $\mathbb{E}\left[(z_{t+1} - \bar{\delta}_t)^2  \right] $ and its value at the previous iterate. We can substitute $t+1$ by $t$ in \eqref{eq:per_iterate_td_average3} to write
\begin{align} \label{eq:per_iterate_td_average_previous}
\mathbb{E}\left[(z_{t} - \bar{\delta}_{t-1})^2  \right]& \leq (1-\beta)\E{(z_{t-1} - \bar{\delta}_{t-2})^2}+ \frac{2 L_V}{\beta} \\
&\quad \times \left[ \alpha^2 (G_{\delta}^2 G_{V}^2 +\lambda^2K^2) + \epsilon^2\right]\nonumber
 + 2 \beta^2 \sigma_{\delta}^2 \; ,
\end{align}
Substituting \eqref{eq:per_iterate_td_average_previous} into the right-hand side of \eqref{eq:per_iterate_td_average3} yields
\begin{align} \label{eq:per_iterate_td_average4}
\!\!\!\!\!\mathbb{E}\!\left[\!(z_{t+1} \!- \! \bar{\delta}_t)^2  \right] &  \!\leq \!(1\!-\!\beta)^2\E{\!(z_{t-1}\! -\! \bar{\delta}_{t-2})^2} \!+ \![1 \!+\! (1\!-\!\beta)]\!\!\! \\
& \quad \times \Big\{ \!\frac{2 L_V}{\beta}\! \! \left[ \alpha^2 (G_{\delta}^2 G_{V}^2 \! +\!\lambda^2K^2) \!+ \! \epsilon^2\right]\!\! + \!2 \beta^2 \sigma_{\delta}^2\! \Big\}.\nonumber
\end{align}
We can recursively apply the previous two steps backwards in time to the initialization to obtain
\begin{align} \label{eq:per_iterate_td_average5}
\mathbb{E}\left[(z_{t+1} - \bar{\delta}_t)^2  \right]  &\leq (1-\beta)^{t+1}{\!(z_{0}\! -\! \bar{\delta}_{-1})^2} +  \!\!\sum_{u=0}^t\!(1\!-\!\beta)^u\!\Big\{\!\frac{2 L_V}{\beta}\! \!  \nonumber\\
&\quad\times \!\! \left[ \alpha^2 \!(G_{\delta}^2 G_{V}^2 \!\!+\!\lambda^2\!K^2) \!+ \! \epsilon^2\right]\!\!+\! 2 \beta^2 \!\sigma_{\delta}^2 \!\Big\} ,\!
\end{align}
In \eqref{eq:per_iterate_td_average5}, the first term on the left-hand side vanishes due to the initialization $z_0=0$ and the convention $\delta_{-1} = 0$. Moreover, the finite geometric sum may be evaluated, provided $\beta<1$, as $\sum_{u=0}^t(1-\beta)^u = [1 - (1-\beta)^t]/\beta$. The numerator in this simplification is strictly less than unit, which means that the right-hand side of \eqref{eq:per_iterate_td_average5} simplifies to 
\begin{align} \label{eq:per_iterate_td_average_final}
\mathbb{E}\left[(z_{t+1} - \bar{\delta}_t)^2  \right] 
 &\leq \frac{2 L_V}{\beta^2}\! \! \left[ \alpha^2 (G_{\delta}^2 G_{V}^2 +\lambda^2K^2) \!+ \! \epsilon^2\right] \!\!+ \!2 \beta \sigma_{\delta}^2 \nonumber\\
& = \ccalO\left(\frac{\alpha^2 + \eps^2}{\beta^2} + \beta\right)
\end{align}
With this relationship established for the auxiliary sequence $z_t$, we shift gears to addressing the evolution of the value function sub-optimality $\|V_t - V^*\|_{\ccalH}$ in expectation. Begin by using the fact that the Hilbert-norm regularizer $(\lambda/2)\|V\|_{\ccalH}^2$ in \eqref{eq:main_prob} implies the objective $J(V)$ is strongly convex, i.e.
\begin{equation}\label{eq:strong_cvx}
 \frac{\lambda}{2} \|V_t - V^* \|_{\ccalH}^2 \leq J(V_t) - V(V^*) \; ,
\end{equation}
together with the expression in Lemma \ref{lemma:iterate_relations}\eqref{lemma:value_function_suboptimality} regarding the value function sub-optimality, assuming constant learning rates and compression budget, i.e. $\alpha_t = \alpha, \beta_t=\beta, \epsilon_t = \epsilon$, to write
\begin{align}\label{eq:value_function_suboptimality_strong_cvx}
&\mathbb{E}\left[\|V_{t+1} - V^* \|_{\ccalH}^2 \given \ccalF_t \right] \nonumber\\
&\quad \leq \left(1 +  \frac{\alpha^2}{\beta} G_V^2 - \alpha \lambda \right) \| V_t - V^* \|_{\ccalH}^2 + 2 \eps \|V_t -V^* \|_{\ccalH} \nonumber\\
&\qquad+ \alpha^2  \sigma_V^2  +  \beta  \mathbb{E}\left[(z_{t+1} - \bar{\delta}_t  )^2 \given \ccalF_t \right] \; .
\end{align}
Consider the total expectation of \eqref{eq:value_function_suboptimality_strong_cvx} with choice of compression budget $\eps=C\alpha^2$ for some arbitrary constant $C>0$, the fact that  $\|V_t - V^* \|_{\ccalH} \leq 2 K$,  apply \eqref{eq:per_iterate_td_average_final} to the last term on the right-hand side of \eqref{eq:value_function_suboptimality_strong_cvx}, and substitute in regularizer $\lambda =   G_V^2 \alpha/\beta + \lambda_0$ for $\lambda_0 < 1$ to obtain:
\begin{align}\label{eq:mean_value_function_suboptimality}
\mathbb{E}\left[\|V_{t+1} - V^* \|_{\ccalH}^2  \right] 
 & \leq \left(1 -\lambda_0 \right)\E{ \| V_t - V^* \|_{\ccalH}^2}  \nonumber\\
&\qquad + \alpha^2 ( \sigma_V^2 + 4 C K) + 2 \beta^2 \sigma_{\delta}^2 \\
&\qquad  + \!   \frac{2 L_V}{\beta}\! \! \left[ \alpha^2 (G_{\delta}^2 G_{V}^2\! \!+\!\lambda^2K^2) \!+ \!C^2 \alpha^4\right] \! .\nonumber
\end{align}
{To establish that  $\liminf \| V_t - V^* \|_{\ccalH}^2 $ is a finite constant determined by $\lambda_0$ and the constant terms on the right-hand side of \eqref{eq:mean_value_function_suboptimality}, which we define as $R:= \alpha^2 ( \sigma_V^2 + 4 C K) + 2 \beta^2 \sigma_{\delta}^2   +   \frac{2 L_V}{\beta}\! \! \left[ \alpha^2 (G_{\delta}^2 G_{V}^2\! \!+\!\lambda^2K^2) \!+ \!C^2 \alpha^4\right]$, suppose that it is not, i.e., that the following holds true:
\begin{align}\label{eq:contradiction_hypothesis}
\liminf_t \mathbb{E}\left[\|V_{t} - V^* \|_{\ccalH}^2  \right] > \frac{R}{\lambda_0}
\end{align}
Then there exists some time index $t_0<\infty$ and some $\delta>0$ such that 
\begin{align}\label{eq:contradiction_hypothesis2}
 \mathbb{E}\left[\|V_{t} - V^* \|_{\ccalH}^2  \right] > \frac{R}{\lambda_0} + \delta
\end{align}
for all $t \geq t_0$. Note that \eqref{eq:contradiction_hypothesis2} may be rearranged to equivalently be stated as
\begin{align}\label{eq:contradiction_hypothesis_rearranged}
\lambda_0\mathbb{E}\left[\|V_{t} - V^* \|_{\ccalH}^2  \right] - \lambda_0 \delta > R 
\end{align}
Let's substitute upper-bound for $R$ stated in \eqref{eq:contradiction_hypothesis_rearranged} into \eqref{eq:mean_value_function_suboptimality}:
\begin{align}\label{eq:mean_value_function_suboptimality2}
\mathbb{E}\left[\|V_{t+1} - V^* \|_{\ccalH}^2  \right] 
 & \leq \left(1 -\lambda_0 \right)\E{ \| V_t - V^* \|_{\ccalH}^2}  + R \nonumber \\
 & < \E{ \| V_t - V^* \|_{\ccalH}^2}  - \lambda_0 \delta \nonumber \\
 & \leq \E{ \| V_t - V^* \|_{\ccalH}^2}
\end{align}
where we have cancelled a common factor of $\lambda_0\mathbb{E}\left[\|V_{t} - V^* \|_{\ccalH}^2  \right]$ from the right-hand side, and upper-estimated $  - \lambda_0 \delta$ by null. Therefore, under the hypothesis that $\liminf_t \mathbb{E}\left[\|V_{t} - V^* \|_{\ccalH}^2  \right] > {R} / {\lambda_0}$, by \eqref{eq:mean_value_function_suboptimality2}, $\mathbb{E}\left[\|V_{t} - V^* \|_{\ccalH}^2\right]$ decreases monotonically to null. This is a contradiction. Therefore, we must have that the hypothesis \eqref{eq:contradiction_hypothesis} is false, and hence 
\begin{align}\label{eq:mean_convergence_subsequence_final}
&\liminf_{t\rightarrow \infty} \E{\|V_{t}\! -\! V^*\|_{\ccalH}^2} =R\nonumber\\
&\quad =\ccalO\left(\alpha^2 + \beta^2 + \frac{\alpha^2}{\beta}\left[1 + \alpha^2 + \frac{\alpha}{\beta} + \frac{\alpha^2}{\beta^2}\right] \right) \; .
\end{align}}
When $\alpha=\beta$, the posynomial of the learning rates on the right-hand side of \eqref{eq:mean_convergence_subsequence_final} simplifies to be $\ccalO(\alpha + \alpha^2+\alpha^3){=\ccalO(\alpha)} \ ${ for $\alpha \in (0,1)$} as stated in \eqref{eq:constant_stepsize_convergence} (Theorem \ref{theorem:constant_stepsize_convergence}).
\subsection{Proof of Corollary \ref{corollary1}}\label{apx_corollary1}

We prove Corollary \ref{corollary1}: In Theorem 3 of \cite{koppel2019parsimonious}[Appendix D.1], it is established for a nonparametric stochastic program without any compositional structure that the effect of sparse subspace projections on the functional stochastic gradient sequence in an RKHS is to yield a function sequence of finite model order, provided a constant algorithm step-size and compression budget are used. The proof of Corollary \ref{corollary1} is nearly identical: the same projection operator is used and the same compactness properties of the state and action spaces apply. The only point of departure is that a distinct deterministic bound is needed on the functional stochastic quasi-gradient for all $\{\bbx_t, \pi(\bbx_t), \bby_t\}$, i.e., to apply the reasoning following equations (74) in \cite{koppel2019parsimonious}[Appendix D.1], we require the existence of a deterministic constant $D$ such that $|[\gamma\kappa(\bby_t,\cdot) - \kappa(\bbx_t, \cdot) ] z_{t+1} | \leq D$ for all $\{\bbx_t, \pi(\bbx_t), \bby_t\}$. We establish such an upper-estimate. To do so, we first establish that the auxiliary sequence $z_t$ stated in \eqref{eq:td_average} is bounded, i.e.








%
%
%
%

\begin{proposition}\label{prop2}
 The auxiliary sequence $z_t$ [cf. \eqref{eq:td_average}] is upper-bounded when used with constant step-size $\beta_t = \beta$:
\begin{equation}\label{eq:prop2}
|z_{t}| = (\gamma+1)K + R_{\max} \text{ for all } t
\end{equation}
\end{proposition}
\begin{myproof}
We pursue a proof by induction. First, the base case: with $V_0 = 0$, we have $|z_1| \leq \beta R_{\max} \leq (\gamma+1)K +  R_{\max}$ making use of the bound on $V_t$ for all $t$ in \eqref{eq:bounded_iterates} and the fact that the step-size is less than unit. 
Now, the induction step: assume the prior bound holds for $z_u$ for $u\leq t$. Write for $z_{t+1}$
\begin{equation}\label{eq:induction_zt}
|z_{t+1}| = (1-\beta) |z_t | + \beta | \delta_t| \leq  (\gamma+1)K + R_{\max}
\end{equation}
where in the last inequality we apply the induction hypothesis together with the upper-estimate on the temporal difference $\delta_t \leq (\gamma+1)K + R_{\max}$.
\end{myproof}

By making use of Proposition \ref{prop2} together with the bound on the reproducing kernel map (Assumption \ref{as:first}), we have the following uniform deterministic bound:

\begin{align}\label{eq:quasi_gradient_bound}
|[\gamma\kappa(\bby_t,\cdot) \!-\! \kappa(\bbx_t, \cdot) ] z_{t+1} |  &\leq X (\gamma \!+\!1) [(\gamma+1) K\! + \!R_{\max} ]  \nonumber \\
&:=D \text{ for all } \{\!\bbx_t, \pi(\!\bbx_t), \bby_t\!\}
\end{align}

Then, we may apply the same reasoning as that of Appendix D.1 of \cite{koppel2019parsimonious} to conclude that the number of Euclidean balls of radius $d=\eps/D$ needed to cover the space $\phi(\ccalX) = \kappa(\ccalX, \cdot)$ is finite, where $\eps$ is a constant as in \eqref{eq:constant_stepsize_parameters}. See \cite{1315946} for further details. Therefore, for Algorithm \ref{alg:pkgtd}, there exists a finite $M^{\infty} < \infty$ such that the model order $M_t \leq M^{\infty}$ for all $t$.



\end{document}